\documentclass[12pt,a4paper]{article}
\usepackage[utf8]{inputenc}
\usepackage[T1]{fontenc}
\usepackage{amsmath,amssymb,amsthm}
\usepackage[margin=2.6cm]{geometry}
\usepackage{tikz}
\usetikzlibrary{arrows.meta,positioning,fit,shapes.geometric,backgrounds}
\usepackage{natbib}
\usepackage{graphicx}
\usepackage{xcolor}
\usepackage{hyperref}
\hypersetup{colorlinks=true,linkcolor=blue!55!black,citecolor=blue!55!black}

\newtheorem{theorem}{Theorem}[section]
\newtheorem{proposition}[theorem]{Proposition}
\newtheorem{corollary}[theorem]{Corollary}
\newtheorem{lemma}[theorem]{Lemma}
\theoremstyle{definition}
\newtheorem{definition}[theorem]{Definition}
\newtheorem{assumption}[theorem]{Assumption}
\theoremstyle{remark}
\newtheorem{remark}[theorem]{Remark}

\newcommand{\Nodes}{\mathcal{N}}
\newcommand{\Arcs}{\mathcal{A}}
\newcommand{\Grph}{\mathcal{G}}
\newcommand{\Vv}{\mathcal{V}}
\newcommand{\Ee}{\mathcal{E}}
\newcommand{\indep}{\perp\!\!\!\perp}
\newcommand{\Below}{\mathrm{Below}}
\newcommand{\Out}{\mathrm{Out}}
\newcommand{\ain}{a_{\mathrm{in}}}
\newcommand{\tw}{\operatorname{tw}}
\newcommand{\Tu}{T^{\mathrm{u}}}       
\newcommand{\Lg}{\mathrm{L}}

\title{Graphical Models Induced by Irreversible Event Accumulation on Trees}
\author{Rodrigo Labouriau\\[2pt]
 {\small Department of Mathematics, Aarhus University}
\footnote{e-mail: rodrigo.labouriau@math.au.dk and rodrigo.labouriau@rlstatlab.com}
}
\date{October 2026}

\begin{document}
\maketitle

\begin{abstract}

Decomposable graphical models are usually adopted for their computational convenience rather than derived from a data-generating mechanism. We show that a special class of them arises naturally from the irreversible accumulation of events along a branching history. Each event is attached to an edge of a latent rooted tree and is inherited by everything below that edge; an observation records which events have occurred at a randomly chosen node, not necessarily a leaf. When the tree is a path, this is the cumulative coding of an ordinal variable, whose indicators form a Markov chain; in the general tree case the Markov chain is replaced by a junction tree.

The conditional independence graph of the event indicators is always decomposable and, in fact, is the line graph of the tree (maximal cliques correspond to nodes where two or more edges meet and separators are the single events on the edges joining two such nodes). The graph determines the tree up to the choice of root. A connected graph arises in this way if, and only if, it is chordal, any two of its maximal cliques share at most one vertex, no vertex is in three maximal cliques, and the graph is such that any two maximal cliques share at most one vertex. This provides a simple graphical test for compatibility with the mechanism (necessary but not sufficient, as irreversibility also constrains the distribution's support).

The structure breaks down if we look at just the leaves, and measurement error destroys the exact screening and leads to a denser graph rather than simply a weaker dependence. We discuss inference, provide diagnostics for departures such as recurrent events, illustrate them on some simulated data and discuss their application in the context of mutation accumulation and lineage tracing where both ancestral and descendant types are sampled together.

\end{abstract}

\medskip
\noindent\textbf{Keywords:} decomposable graphical model; junction tree; event
accumulation; irreversibility; conditional independence; perfect phylogeny;
cancer progression; lineage tracing.

\newpage
\tableofcontents
\newpage

\section{Introduction}
\label{sec:intro}

Multivariate binary data generated by the accumulation of events along
historical trajectories arise in a wide range of scientific areas.  Classic
examples include the presence or absence of mutations in molecular evolution,
somatic mutations in tumour progression, and the transmission of cultural or
linguistic traits.  In such systems, events typically occur irreversibly and are
inherited by everything that descends from the point at which they occurred,
producing strong and highly structured dependence patterns among the collected variables.

A common approach to modelling dependence in multivariate binary data is
through undirected graphical models, and within this system decomposable (or
chordal) graphical models play a central role.  Their importance stems largely
from computational factors: likelihoods factorise according to cliques
and separators, statistical inference becomes tractable, and model-selection
procedures can be implemented efficiently \citep{Lauritzen1996, Edwards2000,
Abreu2010}.  For these reasons decomposable models are often adopted as
convenient approximations or as structural restrictions. However, considerably less
attention has been paid to the question of when decomposability should arise
\emph{naturally} from an underlying data-generating mechanism.

The present work considers this question in the context of irreversible event
accumulation along tree-structured trajectories.  We consider a general
probabilistic model in which each event may occur at most once on a latent
rooted tree and, once realised, is retained by every trajectory that continues
past it.  Each observation is a single \emph{lineage}, sampled at an arbitrary
stage of its progression, and records the presence or absence of the
accumulated events.  This abstract formulation incorporates essential features of a
broad class of cumulative historical processes while remaining independent of
any given application domain.

\paragraph{The model in standard statistical terms.}
Because the vocabulary of the applications may be unfamiliar, we state at once
what the model is when the vocabulary is removed.  A lineage is a directed path
in the tree beginning at the root (i.e., the history of one unit) and, since such a
path is determined by its endpoint, \emph{a lineage is identified with its endpoint node}.
The model therefore draws a single node $Z$ of a tree and reports, for each of the $p$
arcs, whether that arc lies on the root-to-$Z$ path.  Equivalently, it reports
the indicator vector of a \emph{laminar} (nested) family of events, of which
the tree is nothing but the Hasse diagram.  When the tree is a path, this is
exactly the cumulative (``thermometer'') coding of an ordinal variable, and the
main theorem reduces to the familiar fact that cumulative indicators form a
Markov chain; a general tree replaces the \emph{ordinal} variable by a
\emph{hierarchical} one and the Markov chain by a junction tree.  The
observation scheme is likewise standard: each unit is inspected once, at a
random point in its progression, which is cross-sectional (current-status)
sampling of an irreversible progressive process, here with a branching rather
than a linear state space.  Section~\ref{sec:terminology} makes these
statements precise; no knowledge of the applications is required at any point
before Section~\ref{sec:lineage}.

Our main result (Theorem~\ref{thm:main}) characterises the conditional
independence structure induced by this model. The Markov graph of the observed
variables is necessarily decomposable and exhibits a highly constrained clique
structure: maximal cliques intersect in at most one variable, and each variable
belongs to at most two maximal cliques. In particular, the graph admits a
junction-tree representation with singleton separators. Decomposability is
therefore not simply a convenient modelling assumption in this setting, but a
structural consequence of irreversible event accumulation on a tree.

Beyond the structural theorem, these results establish a direct correspondence
between the topology of the event tree and the structure of the induced
graphical model. Maximal cliques correspond to the nodes of the event tree at
which two or more edges meet, while the arc-events joining two such nodes act as
separator variables linking adjacent clique blocks. Thus, the resulting junction
tree reproduces the event tree with its nodes of degree one removed; when every
internal node branches, this is exactly the branching skeleton of the
accumulation process. We complement this with a \emph{converse}
(Theorem~\ref{thm:realis}): the three structural properties of
Theorem~\ref{thm:main} are not only necessary but also sufficient. They
therefore characterise exactly the class of graphs realisable by accumulation
trees, and turn these conditions into a decidable criterion of graphical
compatibility.

\paragraph{A point of modelling emphasis.}
It is tempting to model each observation as a \emph{complete} lineage, that is,
a root-to-leaf path.  We show in Section~\ref{sec:model} that
this apparently natural choice is in fact \emph{degenerate}: complete lineages
satisfy exact flow constraints which collapse the induced graph and
destroy the clique structure described above.
Proposition~\ref{prop:leafcollapse} makes the collapse explicit: under
leaf-only sampling, the Markov graph reduces to a disjoint union of cliques
formed by sibling arcs pointing at leaves, and every other variable becomes an
isolated vertex.  Allowing each lineage to be observed at an arbitrary stage of
progression, so that intermediate states carry positive probability, removes
the degeneracy and yields the structure we describe.  It is also the more
faithful description of applications such as mitochondrial lineage tracing, in
which the sample may contain lineages caught at intermediate stages, not only
at their endpoints \citep{Amorim2019, Ludwig2019}. See Section~\ref{sec:lineage}.

These results suggest a theoretical explanation for the empirical success of
tree-based and decomposable graphical models in systems driven by cumulative
events.  They also suggest natural statistical diagnostics: deviations from the
structural properties derived here may indicate mechanisms such as
recombination, recurrent events, mixtures of populations, or observation noise
\citep{Didelot2010}.

\paragraph{Related work.}
The model sits at the confluence of a number of established frameworks, and stating
the relationship sharpens the contribution.  When each event occurs at most once
and is inherited by all descendants, the units-by-events matrix is a
\emph{perfect phylogeny} in the sense of the infinite-sites model
\citep{Gusfield1991}: the events are binary characters gained once and never
lost, and their compatibility is what forces a tree.  The same ``accumulate,
never lose'' assumption underlies a large family of tumour-progression models:
\emph{oncogenetic trees} \citep{Desper1999, SzaboBoucher2002} restrict
accumulation to a rooted tree; \emph{conjunctive Bayesian networks}
\citep{BeerenwinkelErikssonSturmfels2007, BeerenwinkelSullivant2009, Gerstung2009}
generalise this to a partial order; and \emph{mutual hazard networks}
\citep{Schill2020} replace hard constraints by promoting and inhibiting pairwise
rates.  Direct cross-sectional inference of such progression models is provided
by CAPRESE and CAPRI \citep{OldeLoohuis2014, Ramazzotti2015} (the former for tree
models, the latter for more general progression DAGs) and the family is
surveyed, with software, by \citet{DiazUriarte2022}; a complementary line infers
recurrent trajectories and exclusiveness patterns from cohorts of already-inferred
tumour mutation trees \citep{Luo2023}.  These methods are concerned chiefly with
\emph{inferring} the order of events and their waiting-time dynamics from
cross-sectional data (or from inferred trees).  The present
work asks the complementary, structural question, \emph{given} a tree, what
conditional-independence graph does it impose on the events?, and answers it
with a decomposability theorem and its exact converse
(Theorems~\ref{thm:main} and~\ref{thm:realis}).  The graphical signature thereby
becomes a falsifiable, decidable criterion complementing those estimation
methods (Section~\ref{sec:diagnostics}).

\paragraph{Staged trees and chain event graphs.}
Event trees also underlie staged trees and chain event graphs
\citep{SmithAnderson2008, CollazoGorgenSmith2018}, a class of discrete models
for processes that unfold asymmetrically. There, the event tree lists the
possible developments of a unit, the root-to-leaf paths are the atoms of the
sample space, each non-leaf vertex carries a probability distribution over its
outgoing edges, and vertices with identical transition probabilities are grouped
into stages. Conditional independence, typically of a context-specific kind, is
encoded by the stage structure and read off the graph through separation
theorems \citep{SmithAnderson2008}; distinct staged trees defining the same
statistical model have been characterised by \citet{GorgenSmith2018}, and a
decomposable subclass, with rational maximum likelihood estimators, has been
identified by \citet{AlexandrDuarteVill2024}.

The present model fits into this framework, but the question it asks is
different. Attaching an additional ``stop'' edge to every non-leaf node of $T$
yields an event tree whose root-to-leaf paths correspond one-to-one to the nodes
of $T$. Under Assumption~\ref{ass:pos}, the law of $Z$ is then a
probability tree on this augmented tree, with strictly positive transition
probabilities and no stage identifications; this accounts for the $p$ free
parameters found in Section~\ref{sec:dimension}. Without the stop edges, the
probability tree on $T$ itself describes the leaf-only sampling of
Section~\ref{sec:degenerate}. Rather than reading independence from a stage structure,
we fix one particular vector of random variables, the $p$ arc indicators, and
determine its ordinary undirected Markov graph. Theorem~\ref{thm:main} shows that
this graph depends only on the topology of $T$ and not on the transition
probabilities, provided they are positive; in particular, imposing any stage
structure on the augmented tree leaves it unchanged. This is special to the
accumulation setting: for a general event tree the topology does not by itself
determine conditional independence, which is carried by the stage structure.
Similarly, where
\citet{GorgenSmith2018} ask which trees define the same model,
Proposition~\ref{prop:recovery} and Remark~\ref{rem:rooting} answer the coarser
question of which trees define the same Markov graph: exactly those sharing the
same undirected tree with the same assignment of events to edges. Finally,
decomposability is not built into the model class here, as it is in
\citet{AlexandrDuarteVill2024}; it follows from the accumulation mechanism for
every positive distribution on the nodes.

\paragraph{Organisation.}
Section~\ref{sec:model} introduces the event-accumulation model and the induced
Markov graph.  Section~\ref{sec:main} states the main structural theorem.
Section~\ref{sec:hinge} establishes the hinge-screening lemma and
Section~\ref{sec:proof} gives the direct proof.  Section~\ref{sec:construction}
describes the precise construction of $\Grph(T)$.
Section~\ref{sec:identif} treats identifiability and proves the realisability
converse.  Section~\ref{sec:inference} discusses consequences for inference and
structure learning.  Section~\ref{sec:lineage} instantiates the framework in
lineage tracing.  Section~\ref{sec:diagnostics} develops diagnostics and
discusses departures from the model, and Section~\ref{sec:discussion}
concludes.

\section{The event-accumulation model}
\label{sec:model}

\subsection{The event tree}

Let $T=(\Nodes,\Arcs)$ be a finite rooted directed tree, where $\Nodes$ denotes
the set of nodes and $\Arcs\subseteq\Nodes\times\Nodes$ the set of directed
edges, oriented away from the root.  Following the usual convention for
directed graphs we call the elements of $\Arcs$ \emph{arcs} rather than edges,
reserving ``edge'' for the undirected graphs of Section~\ref{sec:main}; an arc
$a=(u,v)$ has \emph{tail} $u$ and \emph{head} $v$, and we say that $a$
\emph{leaves} $u$ and \emph{enters} $v$.  The \emph{out-degree} $d(v)$ of a
node is the number of arcs leaving it; every non-root node is entered by
exactly one arc, written $\ain(v)$.  For $u,v\in\Nodes$ write
\[
  u \preceq v
\]
if $u$ lies on the unique directed path from the root to $v$; this is the
\emph{ancestral partial order} on $\Nodes$.  For a node $v$, the subtree rooted
at $v$ is
\[
  S(v) \;=\; \{\, z\in\Nodes \,:\, v\preceq z \,\}.
\]
A node is a \emph{leaf} if $d(v)=0$ and \emph{internal} otherwise; it is a
\emph{branching node} if it is the root or has $d(v)\ge 2$, that is, if it is a
point at which the history either starts or divides.

Let $\{E_1,\dots,E_p\}$ be the event types.  Each event $E_i$ is associated
with a unique arc
\[
  a_i=(u_i,v_i)\in\Arcs ,
\]
and distinct events are attached to distinct arcs.  We take $i\mapsto a_i$ to be
a \emph{bijection} between $\{1,\dots,p\}$ and $\Arcs$, so that every arc carries
exactly one observed event and $p=|\Arcs|$ (the \emph{one-event-per-arc}
convention; it is what yields $|\Nodes|=p+1$ in Remark~\ref{rem:bijection} and
the separation arguments below, in which every intervening arc must carry an
indicator).  Forgetting the orientation of the arcs gives the \emph{underlying
undirected tree} $\Tu$, of which each event indexes an edge.
We write $S_i:=S(v_i)$ for the subtree below the head of $a_i$.  Since the arc
$a_i$ is the only way into $v_i$, any history that reaches $v_i$ has
necessarily traversed $a_i$ and hence realised $E_i$; and having realised
$E_i$, it cannot subsequently un-realise it.  These two facts are what we mean
by \emph{persistence} and \emph{irreversibility}: an event, once realised, is
retained for the remainder of the history and is inherited by everything below
the point at which it occurred.  In the model there is no ``forgetting'' step
and no arc that removes an event; formally, irreversibility is the statement
$Y_i=\mathbf 1[Z\in S_i]$ of \eqref{eq:Y} below, that the indicator of $E_i$ is
the indicator of a \emph{subtree}, rather than of an arbitrary set of nodes.

\subsection{Terminology: what a lineage is}
\label{sec:terminology}

The word \emph{lineage} is used throughout in exactly one sense, which we now
fix.

\begin{definition}[Lineage]
\label{def:lineage}
A \emph{lineage} is a directed path in $T$ beginning at the root.  It is the
history of a single unit: the ordered list of events that unit has undergone,
in the order in which it underwent them.  A lineage is \emph{complete} if it
ends at a leaf; otherwise it is a history still in progress.
\end{definition}

Three consequences of Definition~\ref{def:lineage} orient the reader who has no
background in the applications.

\emph{(i) A lineage is identified with its endpoint node.}  A directed path
beginning at the root is determined by its endpoint.  Specifying a lineage is therefore the same as
specifying the node $z\in\Nodes$ it has reached, and ``observing a lineage at
an arbitrary stage of its progression'' means nothing more than drawing a node.
The whole probabilistic content of Section~\ref{sec:obs} is a distribution on
$\Nodes$.  See Figure~\ref{fig:terminology}(a).

\emph{(ii) $Y$ is the indicator vector of a laminar family.}  The subtrees
$S_1,\dots,S_p$ form a laminar (nested) family (Lemma~\ref{lem:laminar}), and
$Y_i=\mathbf 1[Z\in S_i]$.  Conversely, $T$ is precisely the Hasse diagram of
$\{\Nodes,S_1,\dots,S_p\}$ ordered by inclusion.  The tree is thus a
\emph{representation} of how the events nest, not an additional assumption, and
the question the paper answers may be stated with no tree at all: \emph{what
conditional independence structure does the indicator vector of a laminar
family of events acquire under a single draw?}  See
Figure~\ref{fig:terminology}(c).

\emph{(iii) A lineage is not a class of descendants.}  Some of the applications
use ``lineage'' loosely for the set of all continuations below a node.  We do
not: that object is the subtree $S(v)$, and we call it one.  A lineage is a
single path, and an observation is a single lineage.

\begin{figure}[tbp]
\centering
\begin{tikzpicture}[
   >={Stealth[length=1.7mm]},
   nd/.style={circle,draw,gray!70,inner sep=0.4pt,minimum size=5.4mm,
              font=\scriptsize,text=black},
   hit/.style={circle,draw,very thick,fill=black!10,inner sep=0.4pt,
               minimum size=5.4mm,font=\scriptsize},
   ar/.style={->,draw=gray!60},
   on/.style={->,draw,very thick},
   hng/.style={->,draw=red!75!black,line width=1.1pt},
   lb/.style={font=\tiny,inner sep=1.2pt},
   pcap/.style={font=\small,align=center},
   vec/.style={font=\scriptsize,align=center},
   sub/.style={rounded corners=2mm,inner sep=2.5pt,draw=none}]

\begin{scope}
  \node[hit] (a1) at (2.0,3.2)  {$z_1$};
  \node[hit] (a2) at (1.1,2.0)  {$z_2$};
  \node[nd]  (a7) at (3.6,2.0)  {$z_7$};
  \node[nd]  (a3) at (0.0,0.9)  {$z_3$};
  \node[nd]  (a4) at (1.05,0.9) {$z_4$};
  \node[hit] (a5) at (2.3,0.9)  {$z_5$};
  \node[nd]  (a6) at (1.75,-0.3){$z_6$};
  \node[nd]  (a8) at (2.95,-0.3){$z_8$};
  \draw[on] (a1)--(a2) node[lb,midway,left=-1pt] {$E_1$};
  \draw[ar] (a1)--(a7) node[lb,midway,right=-1pt] {$E_6$};
  \draw[ar] (a2)--(a3) node[lb,midway,left=-1pt] {$E_2$};
  \draw[ar] (a2)--(a4) node[lb,midway,left=-2pt] {$E_3$};
  \draw[on] (a2)--(a5) node[lb,midway,above=-2pt] {$E_4$};
  \draw[ar] (a5)--(a6) node[lb,midway,left=-2pt] {$E_5$};
  \draw[ar] (a5)--(a8) node[lb,midway,right=-2pt] {$E_7$};
  \node[font=\scriptsize,anchor=west,red!75!black] at (2.6,0.9) {$Z$};
  \node[vec] at (1.8,-1.15) {$Z=z_5$, incomplete\\[1pt]
                             $Y=(1,0,0,1,0,0,0)$};
  \node[pcap] at (1.8,-2.0) {(a) A lineage is a node};
\end{scope}

\begin{scope}[xshift=7.0cm]
  \node[hit] (b1) at (2.0,3.2)  {$z_1$};
  \node[hit] (b2) at (1.1,2.0)  {$z_2$};
  \node[nd]  (b7) at (3.6,2.0)  {$z_7$};
  \node[nd]  (b3) at (0.0,0.9)  {$z_3$};
  \node[nd]  (b4) at (1.05,0.9) {$z_4$};
  \node[hit] (b5) at (2.3,0.9)  {$z_5$};
  \node[hit] (b6) at (1.75,-0.3){$z_6$};
  \node[nd]  (b8) at (2.95,-0.3){$z_8$};
  \draw[on] (b1)--(b2) node[lb,midway,left=-1pt] {$E_1$};
  \draw[ar] (b1)--(b7) node[lb,midway,right=-1pt] {$E_6$};
  \draw[ar] (b2)--(b3) node[lb,midway,left=-1pt] {$E_2$};
  \draw[ar] (b2)--(b4) node[lb,midway,left=-2pt] {$E_3$};
  \draw[on] (b2)--(b5) node[lb,midway,above=-2pt] {$E_4$};
  \draw[on] (b5)--(b6) node[lb,midway,left=-2pt] {$E_5$};
  \draw[ar] (b5)--(b8) node[lb,midway,right=-2pt] {$E_7$};
  \node[font=\scriptsize,anchor=east,red!75!black] at (1.45,-0.3) {$Z$};
  \node[vec] at (1.8,-1.15) {$Z=z_6$, complete\\[1pt]
                             $Y=(1,0,0,1,\underline{1},0,0)$};
  \node[pcap] at (1.8,-2.0) {(b) Irreversibility: $0\!\to\!1$, never back};
\end{scope}

\begin{scope}[yshift=-6.3cm]
  \node[nd] (c1) at (2.0,3.2)  {$z_1$};
  \node[nd] (c2) at (1.1,2.0)  {$z_2$};
  \node[nd] (c7) at (3.6,2.0)  {$z_7$};
  \node[nd] (c3) at (0.0,0.9)  {$z_3$};
  \node[nd] (c4) at (1.05,0.9) {$z_4$};
  \node[nd] (c5) at (2.3,0.9)  {$z_5$};
  \node[nd] (c6) at (1.75,-0.3){$z_6$};
  \node[nd] (c8) at (2.95,-0.3){$z_8$};
  \draw[ar] (c1)--(c2); \draw[ar] (c1)--(c7);
  \draw[ar] (c2)--(c3); \draw[ar] (c2)--(c4); \draw[ar] (c2)--(c5);
  \draw[ar] (c5)--(c6); \draw[ar] (c5)--(c8);
  \begin{scope}[on background layer]
    \node[sub,fill=blue!14,  fit=(c2)(c3)(c4)(c5)(c6)(c8)] {};
    \node[sub,fill=orange!35,fit=(c5)(c6)(c8)] {};
    \node[sub,fill=green!30, fit=(c7)] {};
  \end{scope}
  \node[lb,anchor=west] at (-0.45,2.30) {$S_1=S(z_2)$};
  \node[lb,anchor=west] at ( 3.20,1.05) {$S_4$};
  \node[lb,anchor=west] at ( 4.05,2.00) {$S_6$};
  \node[vec] at (1.8,-1.15) {$S_4\subset S_1$,\quad $S_6\cap S_1=\emptyset$\\[1pt]
                             $Y_i=\mathbf 1[Z\in S_i]$};
  \node[pcap] at (1.8,-2.0) {(c) A laminar family of subtrees};
\end{scope}

\begin{scope}[xshift=7.0cm,yshift=-6.3cm]
  \node[nd] (d1) at (2.0,3.2)  {$z_1$};
  \node[nd] (d2) at (1.1,2.0)  {$z_2$};
  \node[nd] (d7) at (3.6,2.0)  {$z_7$};
  \node[nd] (d3) at (0.0,0.9)  {$z_3$};
  \node[nd] (d4) at (1.05,0.9) {$z_4$};
  \node[nd] (d5) at (2.3,0.9)  {$z_5$};
  \node[nd] (d6) at (1.75,-0.3){$z_6$};
  \node[nd] (d8) at (2.95,-0.3){$z_8$};
  \draw[ar] (d1)--(d2) node[lb,midway,left=-1pt] {$E_1$};
  \draw[ar] (d1)--(d7) node[lb,midway,right=-1pt] {$E_6$};
  \draw[ar] (d2)--(d3) node[lb,midway,left=-1pt] {$E_2$};
  \draw[ar] (d2)--(d4) node[lb,midway,left=-2pt] {$E_3$};
  \draw[hng](d2)--(d5) node[lb,midway,above=-2pt,red!75!black] {$E_4$};
  \draw[ar] (d5)--(d6) node[lb,midway,left=-2pt] {$E_5$};
  \draw[ar] (d5)--(d8) node[lb,midway,right=-2pt] {$E_7$};
  \begin{scope}[on background layer]
    \node[sub,fill=orange!35,fit=(d5)(d6)(d8)] {};
  \end{scope}
  \node[lb,anchor=west,red!75!black] at (3.20,1.05) {$\mathrm{Below}(4)=\{5,7\}$};
  \node[lb,anchor=west] at (3.95,2.55) {$\mathrm{Out}(4)=\{1,2,3,6\}$};
  \node[vec] at (1.8,-1.15) {$E_4$ is the only way into $S(z_5)$\\[1pt]
                             $Y_{\{5,7\}}\indep Y_{\{1,2,3,6\}}\mid Y_4$};
  \node[pcap] at (1.8,-2.0) {(d) The hinge variable $Y_4$};
\end{scope}
\end{tikzpicture}
\caption{The terminology of Sections~\ref{sec:model} and~\ref{sec:hinge} on a
single event tree (the tree of Figure~\ref{fig:example}).
\textbf{(a)} A \emph{lineage} is a directed path from the root; here it ends at
the internal node $Z=z_5$, so it is \emph{incomplete}.  It has traversed the
arcs $E_1$ and $E_4$, and $Y$ records exactly that.  Because its endpoint determines a root-path, the shaded path and the single node $Z$ carry the
same information: \emph{a lineage is identified with its endpoint node}.
\textbf{(b)} The same lineage, allowed to progress one step further to
$Z=z_6$, now \emph{complete}.  The coordinate $Y_5$ switches from $0$ to $1$
and the coordinates already equal to $1$ stay so: this is
\emph{irreversibility}, and it is the only thing the word means here.  Under
Assumption~\ref{ass:pos} both (a) and (b) are possible observations; under
leaf-only sampling (a) is not, and this is precisely what
Proposition~\ref{prop:leafcollapse} exploits.
\textbf{(c)} Each $Y_i$ is the indicator of a \emph{subtree}, $Y_i=\mathbf
1[Z\in S_i]$, this is the formal content of irreversibility, and subtrees are
nested or disjoint, never partially overlapping (Lemma~\ref{lem:laminar}).
Reading (a) off this panel: $z_5\in S_1$ and $z_5\in S_4$ but $z_5\notin S_6$,
so $Y_1=Y_4=1$ and $Y_6=0$.
\textbf{(d)} The node $z_5$ has exactly one way in, the arc $E_4$.  Its
indicator $Y_4$ is therefore the sole point of articulation between the events
inside $S(z_5)$ and those outside it: conditioning on $Y_4$ makes one of the
two sides constant, and dependence can pass only through it.  This is the
\emph{hinge} of Lemma~\ref{lem:hinge}, and the whole of Theorem~\ref{thm:main}
follows from it.}
\label{fig:terminology}
\end{figure}

\begin{remark}[Two familiar special cases]
\label{rem:familiar}
Both the model and the sampling scheme are standard objects in unfamiliar
dress.

\emph{The linear tree is an ordinal variable.}  Let $T$ be a path
$z_0\to z_1\to\cdots\to z_p$ with $a_i=(z_{i-1},z_i)$, and identify $Z$ with
its index $k\in\{0,1,\dots,p\}$.  Then $Y_i=\mathbf 1[Z\in S_i]=\mathbf 1[k\ge i]$,
which is the cumulative, or ``thermometer'', coding of an ordinal variable
with $p+1$ levels.  Theorem~\ref{thm:main} specialises to the familiar fact
that the cumulative indicators of an ordinal variable form a Markov chain, and
$\Grph(T)$ is the path $Y_1-Y_2-\cdots-Y_p$ (Remark~\ref{rem:chains}).  A
general tree replaces the \emph{ordinal} variable by a \emph{hierarchical}
(nested-category) one, and the Markov chain by a junction tree.  This is the
whole of the generalisation.

\emph{Assumption~\ref{ass:pos} is cross-sectional sampling.}  Each unit is
inspected once, at a random point in its progression, and one records which
events have occurred by then.  This is \emph{current-status} data on an
irreversible progressive process, with a branching rather than a linear state
space.  The root-to-leaf model of Section~\ref{sec:degenerate} corresponds
instead to inspecting only those units that have run to completion, a design
under which the intermediate indicators become deterministic functions of the
terminal one (Proposition~\ref{prop:leafcollapse}); clearly this is a degenerate design.
\end{remark}

\subsection{Observation model}
\label{sec:obs}

Each observation is a single lineage, sampled at an arbitrary stage of its
progression.  Formally, let $Z$ be a random node of $T$, the point reached by
the lineage, which by (i) of Section~\ref{sec:terminology} is the same thing as
the lineage itself, and define the binary random vector $Y=(Y_1,\dots,Y_p)$ by
\begin{equation}\label{eq:Y}
  Y_i \;=\; \mathbf 1\bigl[\, a_i \text{ lies on the root-to-}Z\text{ path}\,\bigr]
      \;=\; \mathbf 1[\, v_i\preceq Z \,]
      \;=\; \mathbf 1[\, Z\in S_i \,].
\end{equation}
Thus $Y_i=1$ precisely when the lineage has progressed far enough to realise
$E_i$.  Observed data consist of independent realisations of $Y$.

\begin{assumption}[Positivity over nodes]
\label{ass:pos}
Every node is a possible observation point:
\[
  \rho(z) \;:=\; \Pr(Z=z) \;>\; 0
  \qquad\text{for every } z\in\Nodes .
\]
\end{assumption}

Assumption~\ref{ass:pos} states that a lineage may be observed at any stage,
including intermediate (non-leaf) nodes and the root.  It is the analytic
counterpart of the modelling remark in the Introduction, and it is exactly the
condition that the root-to-leaf model violates.

\begin{remark}[$Y$ and $Z$ carry the same information]
\label{rem:bijection}
The map $z\mapsto y(z)$ defined by \eqref{eq:Y} is injective: the set of arcs
realised by a lineage is precisely the arc set of the root-to-$Z$ path, which
determines $Z$.  Hence $Y$ and $Z$ are in bijection, and the model is a
saturated multinomial on the $|\Nodes|=p+1$ nodes of $T$.  This observation is
used in Section~\ref{sec:dimension}.
\end{remark}

\subsection{Why complete lineages are degenerate}
\label{sec:degenerate}

Suppose instead that $Z$ is supported only on the leaves of $T$, so that every
observed lineage is complete in the sense of Definition~\ref{def:lineage}.
Then for any arc $a=(u,v)$ whose head $v$ is internal, a lineage reaches $v$ if
and only if it subsequently traverses exactly one outgoing arc of $v$, giving
the exact identity
\begin{equation}\label{eq:flow}
  Y_a \;=\; \sum_{b\,:\,b \text{ leaves } v} Y_b .
\end{equation}
Every internal-arc indicator is therefore a deterministic function of the
others and contributes no conditional-independence constraint of its own.  The
consequence is not merely a loss of information but a complete collapse of the
graph, which the following proposition describes exactly.

\begin{proposition}[Collapse under leaf-only sampling]
\label{prop:leafcollapse}
Suppose $Z$ is supported on the leaves of $T$, with every root-to-leaf path
carrying positive probability, and let $G_{\mathrm{leaf}}$ denote the resulting
Markov graph (Definition~\ref{def:markovgraph}).  Then, for $i\neq j$, the
vertices $i$ and $j$ are adjacent in $G_{\mathrm{leaf}}$ \emph{if and only if}
$a_i$ and $a_j$ are sibling arcs whose heads are both leaves.  Consequently
$G_{\mathrm{leaf}}$ is a disjoint union of cliques,
\[
  G_{\mathrm{leaf}}
  \;=\;
  \biguplus_{v\,:\,\lambda(v)\ge 2}
  K_{\,\{\,a\,:\,a \text{ leaves } v,\ \mathrm{head}(a) \text{ a leaf}\,\}}
  \;\;\uplus\;\;
  \{\text{isolated vertices}\},
\]
where $\lambda(v)$ is the number of leaf children of $v$.  In particular
$G_{\mathrm{leaf}}$ is disconnected whenever $T$ has more than one branching
node, all its separators are \emph{empty}, and it carries none of the clique
structure of Theorem~\ref{thm:main}.
\end{proposition}

\begin{proof}
Fix $i\neq j$ and a configuration $y_{\Vv\setminus\{i,j\}}$ of positive
probability.  The set of root-to-leaf paths compatible with it is the fibre of
$Y_{\Vv\setminus\{i,j\}}$ over that configuration.  If the fibre is a
singleton, then $(Y_i,Y_j)$ is a.s.\ constant given the remaining variables and
is therefore (trivially) conditionally independent.  Adjacency thus requires
some fibre to contain two distinct paths $\pi\neq\pi'$.

Two distinct root-to-leaf paths agree up to their last common node $w$ and then
diverge, so their symmetric difference, as arc sets, is the union of the two
suffixes below $w$, both non-empty and disjoint.  If $\pi,\pi'$ lie in a common
fibre, this symmetric difference is contained in $\{a_i,a_j\}$, which forces
each suffix to consist of a single arc.  Hence $a_i$ and $a_j$ both leave $w$
and both point at leaves.

Conversely, if $a_i,a_j$ are sibling arcs from a node $w$ with leaf heads, the
two paths root-$\to w\to$ head$(a_i)$ and root-$\to w\to$ head$(a_j)$ agree on
every arc other than $a_i,a_j$ and both have positive probability.  On that
fibre $Y_i=1-Y_j$, so $Y_i\not\indep Y_j\mid Y_{\Vv\setminus\{i,j\}}$.
Adjacency of all such sibling pairs at a common $w$ gives the stated cliques,
and every remaining vertex is isolated.
\end{proof}

Applied to the tree of Figure~\ref{fig:example}, Proposition~\ref{prop:leafcollapse}
leaves only the two edges $Y_2\!-\!Y_3$ and $Y_5\!-\!Y_7$, with $Y_1,Y_4,Y_6$
isolated, rather than the ten edges of the three cliques described in
Section~\ref{sec:construction}.  Under Assumption~\ref{ass:pos} the identity
\eqref{eq:flow} fails, a lineage may reach $v$ and stop there, and the
degeneracy disappears.

\subsection{Laminarity}

The following elementary property of subtrees is the only structural fact used
in the proofs.

\begin{lemma}[Laminarity]
\label{lem:laminar}
The family $\{S_1,\dots,S_p\}$ is laminar: for $i\neq j$ exactly one of
\[
  S_i\subseteq S_j,
  \qquad
  S_j\subseteq S_i,
  \qquad
  S_i\cap S_j=\emptyset
\]
holds, according to whether $v_j\preceq v_i$, $v_i\preceq v_j$, or $v_i$ and
$v_j$ are incomparable.
\end{lemma}

\begin{proof}
Subtrees of a rooted tree are nested when their roots are comparable and
disjoint otherwise.  If $v_i\preceq v_j$ then $S_j\subseteq S_i$; if
$v_j\preceq v_i$ then $S_i\subseteq S_j$; and if $v_i,v_j$ are incomparable
then no node lies below both, so $S_i\cap S_j=\emptyset$.
\end{proof}

\subsection{The induced Markov graph}

\begin{definition}[Markov graph induced by $T$]
\label{def:markovgraph}
The \emph{Markov graph induced by $T$} is the undirected graph
$\Grph(T)=(\Vv,\Ee)$ with vertex set $\Vv=\{1,\dots,p\}$, in which distinct
vertices $i,j$ are non-adjacent if and only if
\[
  Y_i \indep Y_j \,\bigm|\, Y_{\Vv\setminus\{i,j\}} .
\]
\end{definition}

Each vertex of $\Grph(T)$ corresponds to a binary variable $Y_i$ indicating the
occurrence of the event $E_i$ attached to an arc of $T$.  Under
Assumption~\ref{ass:pos} every $Y_i$ is non-degenerate, since $S_i$ and its
complement both receive positive mass.

\begin{remark}[Structural zeros and the Markov properties]
\label{rem:zeros}
The law of $Y$ is \emph{never} strictly positive on $\{0,1\}^p$: irreversibility
forces $Y_i=1\Rightarrow Y_j=1$ for every $j$ with $v_j\preceq u_i$, so the
support is a proper subset of the cube.  The classical route from the pairwise
to the global Markov property via the intersection property
(\citealp[Thm.~3.7]{Lauritzen1996}; \citealp{PearlPaz1987}) is therefore
\emph{unavailable in this model}, whatever the observation scheme.  The proof
in Section~\ref{sec:proof} accordingly proceeds in the opposite direction: it
establishes the clique factorisation directly, and reads off the Markov
properties from it.  This is the reason for the structure of the argument, and
it is worth stating explicitly.
\end{remark}

\section{Main theorem}
\label{sec:main}

\begin{theorem}[Structure of the induced Markov graph]
\label{thm:main}
Let $T$ be a finite rooted event tree with $p\ge 2$ arcs and one event per
arc, and let $Y$ be generated as in \eqref{eq:Y} under
Assumption~\ref{ass:pos}.  Then the induced Markov graph $\Grph(T)$ satisfies:
\begin{enumerate}
\item $\Grph(T)$ is chordal (equivalently, decomposable);
\item any two distinct maximal cliques of $\Grph(T)$ intersect in at most one
      vertex;
\item every vertex of $\Grph(T)$ belongs to at most two maximal cliques;
\item $\Grph(T)$ admits a junction-tree representation in which every
      separator is a singleton.
\end{enumerate}
\end{theorem}

(For $p=1$ the graph is a single vertex, its only maximal clique, and the
assertions hold trivially.)  The proof is presented in
Sections~\ref{sec:hinge}--\ref{sec:proof}.  Its main idea is a
single screening identity (Lemma~\ref{lem:hinge}); the remaining steps are the
standard passage from singleton separations to a decomposable factorisation.

\begin{remark}[The Markov graph is a line graph]
\label{rem:linegraph}
Because each event indexes an edge of the underlying undirected tree $\Tu$, the
four clauses admit a compact restatement:
\[
   \Grph(T)\;=\;\Lg(\Tu),
\]
the \emph{line graph} of $\Tu$.  Indeed, two event-variables are adjacent in
$\Grph(T)$ precisely when their arcs share a node of $\Tu$, which is exactly
the identity $\Grph(T)=H$ proved in Section~\ref{sec:proof}.
Conversely, the connected graphs satisfying clauses (1)--(3) are exactly the
line graphs of trees. Indeed, a graph is the line graph of a tree if and only if
it is a connected block graph (every block is complete) in which each cut vertex
lies in exactly two blocks \citep{Harary1963, Chartrand1968}. Block graphs are
precisely the chordal graphs without an induced diamond, $K_4$ minus an edge
\citep{BandeltMulder1986}; clause~(2) is equivalent to having no induced diamond;
and in a block graph the blocks are the maximal cliques, so clause~(3) is the
condition on cut vertices. This is why those three properties reappear as the
realisability criterion of Theorem~\ref{thm:realis}.
\end{remark}

\begin{remark}[Clause (3) is not implied by clause (2)]
\label{rem:threecliques}
The third clause is genuinely independent of the second, and it does the work in
the converse of Section~\ref{sec:identif}.  Consider the graph consisting of
three triangles glued at a common vertex $x$ (Figure~\ref{fig:incompatible}(b)).
It is chordal, its three maximal cliques meet pairwise in exactly $\{x\}$, and
it even admits a junction tree with singleton separators, yet $x$ lies in
\emph{three} maximal cliques, and no accumulation tree produces it: an arc
$(u,v)$ can belong to at most two blocks, namely as an outgoing arc of $u$ and
as the incoming arc of $v$.  Any diagnostic based only on clauses (1), (2), and
(4) would wrongly accept this graph.
\end{remark}

\section{The hinge-screening lemma}
\label{sec:hinge}

The entire theorem rests on one observation, which we state first in words.
Every non-root node $v$ of a tree is entered by exactly one arc.  A lineage can
therefore get below $v$ only by crossing that arc, so the single indicator
attached to it records whether the history has passed into the subtree $S(v)$
or not, and once that binary fact is known, everything on one side of $v$
becomes deterministic.  A variable with this property acts as a \emph{hinge}:
it is the sole point of articulation between the events inside the subtree and
the events outside it, and dependence between the two sides can pass only
through it.  The word is used in this sense throughout, and
Definition~\ref{def:hinge} below fixes it.

\begin{lemma}[Hinge screening]
\label{lem:hinge}
Fix an arc $a_h=(u_h,v_h)$ and write $v=v_h$.  Partition the remaining
indices $\{1,\dots,p\}\setminus\{h\}$ into
\[
  \Below(h) \;=\; \{\, i \,:\, u_i\in S(v) \,\},
  \qquad
  \Out(h)   \;=\; \{1,\dots,p\}\setminus\bigl(\Below(h)\cup\{h\}\bigr).
\]
Then
\[
  Y_{\Below(h)} \;\indep\; Y_{\Out(h)} \,\bigm|\, Y_h .
\]
That is, the single variable $Y_h$ separates the events below $v$ from all
others.
\end{lemma}

\begin{proof}
Observe that $Y_h=\mathbf 1[Z\in S(v)]$.  We condition on the two values of
$Y_h$ and show that in each case one of the two blocks is a constant vector,
which is independent of anything.

\emph{Case $Y_h=0$ (i.e.\ $Z\notin S(v)$).}
For every $i\in\Below(h)$ the tail $u_i$ lies in $S(v)$, hence
$S_i\subseteq S(v)$, and $Z\notin S(v)$ forces $Y_i=0$.  Thus
$Y_{\Below(h)}\equiv 0$ on $\{Y_h=0\}$, and a degenerate vector is
conditionally independent of $Y_{\Out(h)}$.

\emph{Case $Y_h=1$ (i.e.\ $Z\in S(v)$).}
Take any $i\in\Out(h)$.  Since $u_i\notin S(v)$, the arc $a_i$ is not interior
to $S(v)$, so by laminarity (Lemma~\ref{lem:laminar}) either
$S_i\cap S(v)=\emptyset$ or $S_i\supseteq S(v)$.  (The remaining possibility
$S_i\subsetneq S(v)$ is excluded because $v$ has a \emph{unique} incoming arc,
namely $a_h$: any arc whose head lies strictly inside $S(v)$ has its tail in
$S(v)$ as well, hence lies in $\Below(h)$.)  In the first case $Z\in S(v)$
forces $Y_i=0$; in the second case $Z\in S(v)\subseteq S_i$ forces $Y_i=1$.
Either way $Y_i$ is constant on $\{Y_h=1\}$.  Hence $Y_{\Out(h)}$ is a
degenerate vector on this event and is conditionally independent of
$Y_{\Below(h)}$.

In both cases the conditional joint law of
$\bigl(Y_{\Below(h)},Y_{\Out(h)}\bigr)$ given $Y_h$ factorises, which is the
assertion.
\end{proof}

\begin{definition}[Hinge variable]
\label{def:hinge}
For a non-root node $v$, the variable $Y_h$ of Lemma~\ref{lem:hinge} with
$h=\ain(v)$ is the \emph{hinge variable} of $v$: it screens the events strictly
below $v$ from all the others.  Every non-root arc-event is a hinge in this
sense.
\end{definition}

Two remarks on the terminology.  First, screening is universal but being a
junction-tree \emph{separator} is not.  A hinge $Y_h$ is a singleton separator
precisely when its arc joins two nodes of $\Tu$ of degree at least two,
equivalently, when $Y_h$ lies in two maximal cliques
(Corollary~\ref{cor:edges}).  An arc entering a leaf, for instance, screens a
\emph{trivial} subtree ($\Below(h)=\varnothing$) and belongs to a single clique,
so it is a hinge but not a separator; we therefore do not identify all hinges
with separators.  Second, the mechanism arises from a \emph{degeneracy}, not from
averaging. Indeed, conditioning on a hinge does not only weaken the dependence across it; it also
makes one side constant.  This is why the screening is exact, and
also, as Section~\ref{sec:noise} shows, why it is fragile under measurement
error, since a noisy hinge no longer determines either side.

\section{Proof of the theorem}
\label{sec:proof}

Throughout this section, for every \emph{non-leaf} node $v$ of $T$ define the
\emph{block}
\begin{equation}\label{eq:block}
  C(v)
  \;=\;
  \{\, a\in\Arcs : a \text{ leaves } v \,\}
  \;\cup\;
  \begin{cases}
    \{\ain(v)\}, & v\neq\text{root},\\
    \emptyset,   & v=\text{root},
  \end{cases}
\end{equation}
where $\ain(v)$ denotes the unique arc entering $v$.  Let
$\mathcal{B}=\{\,C(v) : v \text{ non-leaf},\ |C(v)|\ge 2\,\}$; the size
condition discards only $C(\text{root})$ in the degenerate case where the root
has out-degree one, in which case that single vertex already lies in the block
of the next node.  Let $H$ be the undirected graph on $\{1,\dots,p\}$ whose
edges are all pairs contained in a common block,
$\Ee(H)=\bigcup_{C\in\mathcal{B}}\binom{C}{2}$.  We prove
Theorem~\ref{thm:main} by showing $\Grph(T)=H$ and that $H$ is decomposable
with singleton separators. The proof is split into four parts.

\paragraph{First: each arc lies in at most two blocks, and the
block-intersection graph is a tree.}
An arc $a=(u,v)$ leaves exactly one node and enters exactly one node.  It
belongs to $C(w)$ as an outgoing arc only for $w=u$, and as the incoming arc
only for $w=v$ (and only if $v$ is non-leaf).  Hence $a$ lies in at most two
blocks, which is clause (3) of Theorem~\ref{thm:main}.

Consequently, two distinct blocks $C(u)$ and $C(w)$ share a vertex if and only
if one of $u,w$ is the parent of the other: an arc lying in both must leave one
and enter the other.  The shared vertex is then the unique arc joining them,
so \emph{every block intersection is a singleton}, which is clause (2).
Moreover, the graph on the blocks, with an edge whenever two blocks meet, is
exactly the subtree of $\Tu$ induced by the nodes of degree at least two: two
blocks $C(u),C(w)$ meet precisely when $u$ and $w$ are adjacent in $\Tu$, sharing
the single arc between them.  Equivalently, it is $T$ with \emph{all} its pendant
vertices removed, including the root when it has out-degree one, whose block
$C(\mathrm{root})$ is then the singleton excluded from $\mathcal B$.  (Removing
only the leaves would wrongly retain such a root.)  Call this tree the
\emph{block tree} and denote it $T^{\circ}$.

\paragraph{Second: $\Grph(T)\subseteq H$ (the factorisation).}
We prove the junction-tree factorisation
\begin{equation}\label{eq:factor}
  \Pr(y)
  \;=\;
  \frac{\prod_{C\in\mathcal B} \Pr\bigl(y_{C}\bigr)}
       {\prod_{s\in\mathcal S} \Pr\bigl(y_{s}\bigr)} ,
\end{equation}
over the blocks, where $\mathcal B$ is the set of blocks (the maximal cliques,
by Step~4) and
$\mathcal S$ is the multiset of the $|\mathcal B|-1$ singleton separators
labelling the edges of the block tree $T^{\circ}$, namely the arcs joining two
nodes of degree at least two.  Only these arcs enter the denominator: an arc
entering a leaf is a private vertex of a single clique and does \emph{not}
separate, so it must not be divided out.  (Ranging the denominator over all
non-root arcs, as a naive reading of ``one separator per non-root node'' would
suggest, over-counts and does not reproduce $\Pr(y)$.)

The proof of \eqref{eq:factor} is by induction on the number of non-leaf nodes
of $T$, simultaneously for all trees with $p\ge2$ arcs and all $\rho$ satisfying
Assumption~\ref{ass:pos}.  If the root is the only non-leaf node, every arc
leaves the root, $\mathcal B=\{C(\mathrm{root})\}$ with $|C(\mathrm{root})|=p\ge
2$, $\mathcal S=\emptyset$, and \eqref{eq:factor} is the identity
$\Pr(y)=\Pr(y_{C(\mathrm{root})})$.  Otherwise choose a non-leaf node
$v\neq\mathrm{root}$ all of whose children are leaves (a deepest non-leaf node),
and let $h=\ain(v)=(u,v)$.  Then $\Below(h)$ consists of the outgoing arcs of
$v$, so $\Below(h)\cup\{h\}=C(v)$, and Lemma~\ref{lem:hinge} gives
\[
  \Pr(y)
  \;=\;
  \frac{\Pr\bigl(y_{C(v)}\bigr)\;
        \Pr\bigl(y_{\Out(h)\cup\{h\}}\bigr)}
       {\Pr(y_h)} ,
\]
where $\Pr(y_h)>0$ for both values of $y_h$ because $\rho(S(v))>0$ and
$\rho(\Nodes\setminus S(v))>0$.  Let $T'$ be the tree obtained from $T$ by
deleting the children of $v$, so that $v$ becomes a leaf, and let
$\rho'(z)=\rho(z)$ for $z\notin S(v)$ and $\rho'(v)=\rho(S(v))$; $\rho'$ satisfies
Assumption~\ref{ass:pos}.  With $Z'=v$ when $Z\in S(v)$ and $Z'=Z$ otherwise,
$Z'\sim\rho'$, and the indicators of the arcs of $T'$ evaluated at $Z'$ are
exactly $Y_{\Out(h)\cup\{h\}}$.  The second factor is therefore the law of the
indicator vector of $(T',\rho')$, a tree with $p'=p-d(v)$ arcs and one non-leaf
node fewer.  If $p'=1$, then $T'$ is the single arc $h$, so $u$ is a root of
out-degree one; the second factor is $\Pr(y_h)$ and cancels, and
$\Pr(y)=\Pr(y_{C(v)})$ is \eqref{eq:factor} for $T$, whose only block is $C(v)$.
If $p'\ge2$, the induction hypothesis factorises the second factor over the
blocks $\mathcal B'$ and separators $\mathcal S'$ of $T'$.  The blocks of $T$
are those of $T'$ together with $C(v)$ (all other blocks are unchanged, and $v$
is a leaf of $T'$), so $\mathcal B=\mathcal B'\cup\{C(v)\}$.  Since $p'\ge2$, the
node $u$ is not a root of out-degree one, hence $|C(u)|\ge2$ and $C(u)\in\mathcal
B'$; the arc $h$ joins the nodes $u$ and $v$, both of degree at least two in
$\Tu$, so it is a separator of $T$, and $\mathcal S=\mathcal S'\cup\{h\}$.
Substituting gives \eqref{eq:factor} for $T$.  Positivity enters this step only
through the denominators; it is essential in Step~3.

A distribution that factorises over the complete subsets of a graph satisfies
the global Markov property with respect to it, positive or not
\citep[Prop.~3.8]{Lauritzen1996}.  Hence every non-edge of $H$ is a non-edge of
$\Grph(T)$, i.e.\ $\Grph(T)\subseteq H$.

\paragraph{Third: $H\subseteq\Grph(T)$ (each block is complete).}
Fix a non-leaf node $v$ and two vertices of $C(v)$; we exhibit a configuration
of the remaining variables, of positive probability, under which the two remain
dependent.  In each case we use Assumption~\ref{ass:pos}, which gives positive
mass to each individual node.

\begin{itemize}
\item \emph{Two outgoing arcs $o,o'$ of $v$.}  Condition on the event that
  every variable except $Y_o,Y_{o'}$ takes the value forcing
  $Z\in\{v,\,v_o,\,v_{o'}\}$ (reach $v$; take no other branch at $v$; do not
  descend below $v_o$ or $v_{o'}$; and follow the unique path to $v$).  Then
  $(Y_o,Y_{o'})\in\{(0,0),(1,0),(0,1)\}$, each with positive probability, while
  $(1,1)$ is impossible.  The conditional law does not factorise, so
  $Y_o\not\indep Y_{o'}\mid Y_{\mathrm{rest}}$.

\item \emph{The incoming arc $h$ and an outgoing arc $o$ of $v$.}  Condition so
  that $Z\in\{u_h,\,v,\,v_o\}$ (reach the tail $u_h$ of $h$; then either stop,
  cross $h$ into $v$, or continue along $o$).  Then
  $(Y_h,Y_o)\in\{(0,0),(1,0),(1,1)\}$, each with positive probability, while
  $(0,1)$ is impossible, so $Y_h\not\indep Y_o\mid Y_{\mathrm{rest}}$.
\end{itemize}

Thus every edge of $H$ is an edge of $\Grph(T)$, and combined with Step 2 we
obtain $\Grph(T)=H$.

\paragraph{Fourth: conclusion.}
By $\Grph(T)=H$ and Step 1, the maximal cliques of $\Grph(T)$ are exactly the
blocks in $\mathcal{B}$ (no block is contained in another, since distinct
blocks meet in at most one vertex, and every block of $\mathcal{B}$ has at least
two vertices); distinct maximal cliques meet in at most one vertex; each vertex
lies in at most two of them; and the block tree $T^{\circ}$ is a junction tree
for $\Grph(T)$ whose separators are the singleton arcs joining two nodes of
degree at least two.  The existence
of a junction tree is equivalent to chordality \citep{Lauritzen1996,
Rose1976}.  This establishes all four assertions of Theorem~\ref{thm:main}.
\hfill$\square$

\begin{remark}[Non-adjacency, made explicit]
\label{rem:nonadj}
The inclusion $\Grph(T)\subseteq H$ can also be read off directly.  If $i,j$
share no block then, in the comparable case, some arc $h$ lies strictly between
$a_i$ and $a_j$ on the root-to-leaf path containing both; in the incomparable
case, some initial branch arc $h$ at their lowest common ancestor lies below
exactly one of them.  In either case $\{h\}$ separates $i$ from $j$ by Lemma~\ref{lem:hinge},
so $i$ and $j$ are non-adjacent.  In particular (and this is the point at
which a natural intuition might mislead) \emph{incomparable sibling arcs leaving a
common node are adjacent}: they lie in the common block $C(v)$ and are coupled
by the single-branch (mutual-exclusion) constraint at $v$.  Incomparability of
arcs is therefore \emph{not} sufficient for conditional independence.
\end{remark}

\section{Construction of the Markov graph}
\label{sec:construction}

The proof already yields an explicit construction of $\Grph(T)$.  We now record
it, and specialise to the case in which every internal node of $T$ branches.

\begin{proposition}[Blocks are the maximal cliques]
\label{prop:blocks}
Under Assumption~\ref{ass:pos}, the maximal cliques of $\Grph(T)$ are
precisely the blocks $C(v)$ of \eqref{eq:block} for non-leaf nodes $v$ with
$|C(v)|\ge 2$.  Moreover, let $v$ be neither the root nor a leaf, with parent
$u$.  If $u$ has degree at least two in $\Tu$ (that is, unless $u$ is a root of
out-degree one), then $\ain(v)$ is a singleton separator: it is the unique
vertex shared between $C(v)$ and $C(u)$.  If $u$ is a root of out-degree one,
$\ain(v)$ lies in the single maximal clique $C(v)$ and is not a separator.
\end{proposition}

\begin{proof}
This is the content of Steps 1, 3 and 4 of the proof of
Theorem~\ref{thm:main}: the blocks are complete (Step 3), no larger set is
complete (Step 2, since any additional vertex is separated from some member of
$C(v)$ by an intervening incoming arc), each pair of blocks meets in at most
the single incoming arc of the deeper node (Step 1), and this arc is the
separator of the junction tree (Step 4).
\end{proof}

\begin{definition}[Branching skeleton]
\label{def:skeleton}
The \emph{branching skeleton} $B(T)$ of $T$ is the tree whose nodes are the
branching nodes of $T$ (the root, and the nodes of out-degree at least two),
each joined to its nearest branching ancestor.  Informally, it is what remains
of $T$ after every stretch along which the history has no choice to make has
been contracted to a single connection: it records \emph{where} the history can
divide, and how those division points are nested, and discards everything else.
\end{definition}

\begin{proposition}[Branching nodes generate maximal cliques]
\label{prop:branching}
Assume in addition that every internal node of $T$ has out-degree at least two
(so that every non-leaf node is a branching node).  Then the maximal cliques of
$\Grph(T)$ are in bijection with the branching nodes of $T$, and the junction
tree of $\Grph(T)$ is isomorphic to the branching skeleton $B(T)$, with edges
labelled by the singleton incoming-arc separators.
\end{proposition}

\begin{proof}
Under the additional hypothesis, the non-leaf nodes are exactly the branching
nodes, $|C(v)|\ge 2$ for all of them, and the block tree $T^{\circ}$ of Step 1
coincides with $B(T)$.  Apply Proposition~\ref{prop:blocks}.
\end{proof}

\paragraph{Structural interpretation.}
The clique structure of $\Grph(T)$ mirrors the branching structure of the event
tree.  Under the branching convention of Proposition~\ref{prop:branching} each
branching node generates a maximal clique ,  in general, every node of degree at
least two does, an out-degree-one internal node contributing a two-vertex clique.
Such a node is a point at which a history may continue in several alternative
ways, and because it may continue in only \emph{one} of them, the corresponding
event indicators are mutually exclusive and hence dependent, this is the
single-branch constraint, and it is what makes the block complete.  The incoming
arc of such a node is its hinge variable (Definition~\ref{def:hinge}): it
separates the dependence structure above from that below.  The junction-tree
separators are then the singleton arc-events joining two nodes of degree at least
two.

\begin{corollary}[Clique size and node degree]
\label{cor:cliquesize}
Let $v$ be a non-leaf node of $T$ with $|C(v)|\ge 2$; equivalently, let $v$ be a
node of degree at least two in $T^{u}$ (every non-leaf node qualifies except a
root of out-degree one). Then $C(v)$ is a maximal clique of $\Grph(T)$, with
cardinality
\[
  |C(v)| \;=\;
  \begin{cases}
    d(v),   & v \text{ the root},\\
    d(v)+1, & v \neq \text{root},
  \end{cases}
  \qquad\text{that is,}\qquad |C(v)| = \deg_{T^{u}}(v).
\]
Conversely, every maximal clique of $\Grph(T)$ equals $C(v)$ for a unique such
node $v$. Under the hypothesis of Proposition~\ref{prop:branching}, these nodes
are exactly the branching nodes of $T$.
\end{corollary}

\begin{proof}
By Proposition~\ref{prop:blocks}, the maximal cliques of $\Grph(T)$ are exactly
the blocks $C(v)$ with $v$ non-leaf and $|C(v)|\ge 2$. By~\eqref{eq:block},
$C(v)$ consists of one variable per outgoing arc of $v$, together with the
incoming arc when $v$ is not the root. Hence $|C(v)| = d(v)$ for the root and
$|C(v)| = d(v)+1$ otherwise, which in both cases is the degree of $v$ in
$T^{u}$. In particular $|C(v)|\ge 2$ fails only for a root of out-degree one.
Uniqueness: distinct blocks meet in at most one vertex (Step~1 of the proof of
Theorem~\ref{thm:main}), whereas each block in $\mathcal{B}$ has at least two
vertices, so $u\neq w$ implies $C(u)\neq C(w)$. Under the hypothesis of
Proposition~\ref{prop:branching}, every non-leaf node, the root included, has
out-degree at least two, so these nodes are exactly the branching nodes.
\end{proof}

\begin{corollary}[Edges, separators, treewidth]
\label{cor:edges}
Let $\Grph(T)$ have maximal cliques $C_1,\dots,C_m$.  Then:
\begin{enumerate}
\item every edge of $\Grph(T)$ lies in \emph{exactly one} maximal clique;
\item exactly $m-1$ vertices lie in two maximal cliques (the separators), and
      the remaining $p-(m-1)$ vertices lie in exactly one; equivalently
      $\sum_{a=1}^{m}|C_a| = p+m-1$;
\item the treewidth of $\Grph(T)$ is
      $\tw\bigl(\Grph(T)\bigr)=\max_a |C_a|-1$, which under the hypotheses of
      Proposition~\ref{prop:branching} equals
      $\max\bigl\{\,d(\mathrm{root})-1,\ \max_{v\neq\mathrm{root}} d(v)\,\bigr\}$.
\end{enumerate}
\end{corollary}

\begin{proof}
(1) If an edge $\{i,j\}$ lay in two distinct maximal cliques $C_a,C_b$, then
$|C_a\cap C_b|\ge 2$, contradicting Theorem~\ref{thm:main}(2); and every edge
lies in at least one block by construction.
(2) The junction tree has $m$ nodes and $m-1$ edges, each labelled by a
distinct separator vertex; distinctness holds because a vertex in two cliques
determines a unique junction-tree edge, and a vertex cannot lie in three
cliques by Theorem~\ref{thm:main}(3).  Counting incidences gives
$\sum_a|C_a| = (p-(m-1))\cdot 1 + (m-1)\cdot 2 = p+m-1$.
(3) For a chordal graph the treewidth is the clique number minus one, and the
clique sizes are given by Corollary~\ref{cor:cliquesize}.
\end{proof}

Part (3) governs the exponential factor in the cost of exact inference on
$\Grph(T)$, and that factor is controlled by the branching degree of the event
tree alone; the total cost grows in addition linearly with the numbers of
cliques and variables (Section~\ref{sec:inference}).

\begin{remark}[Chains of out-degree-one nodes]
\label{rem:chains}
If $T$ contains internal nodes of out-degree one, then under
Assumption~\ref{ass:pos} the arcs along such a chain remain distinct
non-degenerate variables (a lineage may stop between them), and a maximal chain
contributes a path of two-vertex cliques joined by singleton separators.  All
four conclusions of Theorem~\ref{thm:main} are unaffected; only the bijection
between maximal cliques and \emph{branching} nodes requires the branching
hypothesis of Proposition~\ref{prop:branching}.  This is in sharp contrast to
the root-to-leaf model of Proposition~\ref{prop:leafcollapse}, in which chain
arcs are perfectly confounded and collapse to a single variable.
\end{remark}

\begin{figure}[t]
\centering
\begin{tikzpicture}[
   >={Stealth[length=2mm]},
   nd/.style={circle,draw,inner sep=1pt,minimum size=6mm,font=\small},
   yv/.style={circle,draw,inner sep=1pt,minimum size=6.5mm,font=\small},
   every edge/.style={draw,->}]

\begin{scope}
\node[nd] (z1) at (2.1,3.0) {$z_1$};
\node[nd] (z2) at (1.1,1.9) {$z_2$};
\node[nd] (z7) at (3.5,1.9) {$z_7$};
\node[nd] (z3) at (-0.2,0.8) {$z_3$};
\node[nd] (z4) at (1.0,0.8) {$z_4$};
\node[nd] (z5) at (2.3,0.8) {$z_5$};
\node[nd] (z6) at (1.7,-0.4) {$z_6$};
\node[nd] (z8) at (2.9,-0.4) {$z_8$};
\path (z1) edge node[left,font=\scriptsize,pos=.5] {$E_1$} (z2)
      (z1) edge node[right,font=\scriptsize,pos=.5] {$E_6$} (z7)
      (z2) edge node[left,font=\scriptsize,pos=.5] {$E_2$} (z3)
      (z2) edge node[left,font=\scriptsize,pos=.5] {$E_3$} (z4)
      (z2) edge node[right,font=\scriptsize,pos=.5] {$E_4$} (z5)
      (z5) edge node[left,font=\scriptsize,pos=.5] {$E_5$} (z6)
      (z5) edge node[right,font=\scriptsize,pos=.5] {$E_7$} (z8);
\node at (1.7,-1.3) {\small (a) Event tree $T$, one event per arc};
\end{scope}

\begin{scope}[xshift=7.6cm]
\node[yv] (Y6) at (0.0,3.0) {$Y_6$};
\node[yv] (Y1) at (1.5,3.0) {$Y_1$};
\node[yv] (Y2) at (0.5,1.8) {$Y_2$};
\node[yv] (Y3) at (1.7,1.8) {$Y_3$};
\node[yv] (Y4) at (2.9,1.8) {$Y_4$};
\node[yv] (Y5) at (2.3,0.5) {$Y_5$};
\node[yv] (Y7) at (3.7,0.5) {$Y_7$};
\draw (Y6)--(Y1);
\draw (Y1)--(Y2); \draw (Y1)--(Y3); \draw (Y1)--(Y4);
\draw (Y2)--(Y3); \draw (Y2)--(Y4); \draw (Y3)--(Y4);
\draw (Y4)--(Y5); \draw (Y4)--(Y7); \draw (Y5)--(Y7);
\node at (1.9,-1.3) {\small (b) Induced Markov graph $\Grph(T)$};
\end{scope}
\end{tikzpicture}
\caption{The map $T\mapsto\Grph(T)$ in the one-event-per-arc setting.  The
maximal cliques are $C_0=\{Y_1,Y_6\}$ (node $z_1$),
$C_1=\{Y_1,Y_2,Y_3,Y_4\}$ (node $z_2$, out-degree three) and
$C_2=\{Y_4,Y_5,Y_7\}$ (node $z_5$, out-degree two), with singleton separators
$C_0\cap C_1=\{Y_1\}$ and $C_1\cap C_2=\{Y_4\}$.  The junction tree is the path
$C_0-C_1-C_2$, isomorphic to $T$ with its leaves removed.  Under leaf-only
sampling, by contrast, Proposition~\ref{prop:leafcollapse} leaves only the
edges $Y_2\!-\!Y_3$ and $Y_5\!-\!Y_7$.}
\label{fig:example}
\end{figure}

\section{Identifiability and realisability}
\label{sec:identif}
We now ask to what extent the event tree $T$ can be recovered from the induced
graphical model, and, conversely, which graphs arise this way at all.

\subsection{Recovery of the event tree}

\begin{proposition}[Recovery of the unrooted tree]
\label{prop:recovery}
Under Assumption~\ref{ass:pos} and one event per arc, $\Grph(T)$ determines the
underlying undirected tree $\Tu$ up to isomorphism, together with the assignment
of each event to an edge of $\Tu$.  Explicitly, $\Tu$ is reconstructed by
inverting the line-graph map (Remark~\ref{rem:linegraph}): create one tree node
for each maximal clique of $\Grph(T)$; join two clique-nodes by an edge carrying
their common separator variable; and, for every variable lying in a single
clique, attach a fresh pendant node by an edge carrying that variable.
\end{proposition}

\begin{proof}
By Theorem~\ref{thm:main}, the maximal cliques of $\Grph(T)$ are the blocks and
their singleton intersections, the separators, so the junction tree, recoverable
from the graph, is the block tree $T^{\circ}$.  Restoring, for each private
(non-separator) clique vertex, a pendant node joined by an edge carrying that
event returns precisely the degree-one nodes of $\Tu$ deleted in forming
$T^{\circ}$, while separator vertices label the internal edges.  This inverts the
line-graph map of Remark~\ref{rem:linegraph}, and a tree is determined up to
isomorphism by its line graph, so $\Tu$ is recovered.  Each event is placed by
its clique-incidence pattern: two blocks for a separator (an internal edge), one
block for a private event (a pendant edge).
\end{proof}

The observation scheme of Assumption~\ref{ass:pos} improves identifiability
relative to the root-to-leaf model.  Because a lineage may stop at any node,
arcs lying on out-degree-one segments are non-degenerate and appear explicitly
in $\Grph(T)$ (Remark~\ref{rem:chains}); likewise, events occurring before the
first branching are not constant, since the lineage may be observed at or above
the first branching node.  The segments that are non-identifiable under
leaf-only sampling therefore become identifiable here.  Arcs pointing at leaves
of $T$ are recovered as the non-separator vertices of their blocks, so in fact
the whole of $\Tu$ is recovered; the only ambiguity that remains is the rooting,
described next.

\begin{remark}[Rooting]
\label{rem:rooting}
The residual ambiguity is one of \emph{orientation}, and it is more severe than a
choice among blocks.  Because $\Grph(T)=\Lg(\Tu)$ depends on $T$ only through its
underlying undirected tree, \emph{every} rooting of $\Tu$ yields the same Markov
graph.  The root may therefore be any node of $\Tu$, including a degree-one node,
which corresponds to no maximal clique, and orienting all edges away from it
gives a rooted event tree with the same $\Grph(T)$.  Thus $\Grph(T)$ determines
$\Tu$ up to unrooted isomorphism and no more.

Marginal frequencies do not resolve the root either.  On the path $A-B-C$ with
arc-events $e_1=AB$ and $e_2=BC$, rooting at $A$ with node masses
$(\rho_A,\rho_B,\rho_C)=(0.4,0.3,0.3)$ and rooting at $B$ with $(0.6,0.1,0.3)$
both give $\Pr(Y_{e_1}=1)=0.6$ and $\Pr(Y_{e_2}=1)=0.3$, with all node masses
positive; the marginals cannot distinguish them.  What separates them is the
\emph{joint} support: under the first rooting the state $(1,1)$ occurs (the
events are nested), whereas under the second the two events are mutually
exclusive.  Recovering the orientation therefore requires the implication and
mutual-exclusion relations carried by the full distribution, not the graph, and
not the one-way marginals; we do not pursue a formal root-identifiability result
here.
\end{remark}

\subsection{Which graphs arise? An exact characterisation.}

\begin{theorem}[Realisability]
\label{thm:realis}
Let $G$ be a connected undirected graph on $p\ge 2$ vertices.  Then $G$ is
isomorphic to $\Grph(T)$ for some event-accumulation tree $T$ with one event
per arc, observed under Assumption~\ref{ass:pos}, \emph{if and only if}
\begin{enumerate}
\item[(i)]   $G$ is chordal;
\item[(ii)]  any two distinct maximal cliques of $G$ intersect in at most one
vertex;
\item[(iii)] every vertex of $G$ lies in at most two maximal cliques.
\end{enumerate}
\end{theorem}

\begin{proof}
Necessity is Theorem~\ref{thm:main}.  For sufficiency, let $C_1,\dots,C_m$ be
the maximal cliques of $G$ and let $K$ be the \emph{clique-intersection graph}:
vertices $C_1,\dots,C_m$, with an edge whenever $C_a\cap C_b\neq\emptyset$.  By
(ii) each such intersection is a singleton, say ${s_{ab}}$.  By (iii),
distinct edges of $K$ carry distinct separator vertices: if the edges
${C_a,C_b}$ and ${C_a,C_c}$ shared the separator $s$, then $s$ would lie in
the three maximal cliques $C_a,C_b,C_c$.

\emph{$K$ is a tree.}  It is connected because $G$ is.  Since $G$ is chordal, it
possesses a junction tree $J$, necessarily a spanning tree of $K$.  Fix any
separator vertex $s$, lying in exactly the two cliques $C_a,C_b$.  The running
intersection property forces the cliques containing $s$ to induce a connected
subtree of $J$; as there are only two of them, $C_a$ and $C_b$ must be adjacent
in $J$.  Hence every edge of $K$ is an edge of $J$, so $K\subseteq J\subseteq
K$ and $K=J$ is a tree.

\emph{Construction of $T$.}  Root $K$ at an arbitrary clique $C_r$.  For
$a\neq r$ let $s_a$ denote the separator on the edge joining $C_a$ to its
parent in $K$, and set $s_r:=\emptyset$.  Introduce a node $\nu_a$ for each
clique $C_a$.  For each $a$ and each $x\in C_a\setminus{s_a}$:
\begin{itemize}
\item if $x=s_b$ for some child $C_b$ of $C_a$ in $K$, create the arc
$(\nu_a,\nu_b)$ and attach the event $x$ to it;
\item otherwise create a fresh leaf $\ell_x$ and the arc $(\nu_a,\ell_x)$, with
the event $x$ attached to it.
\end{itemize}
The assignment is well defined: $s_b\in C_a$ for every child $C_b$, and
$s_b\neq s_a$ by the distinctness of separators, so the elements of
$C_a\setminus{s_a}$ are used exactly once each.  Every vertex $x$ of $G$
labels exactly one arc: if $x$ lies in two cliques it is some $s_b$ and labels
the single arc $(\nu_{\text{parent}},\nu_b)$; if it lies in one clique $C_a$
then $x\neq s_a$ and it labels $(\nu_a,\ell_x)$.

The resulting $T$ is a rooted tree with one event per arc.  By construction the
node $\nu_a$ has incoming arc $s_a$ (for $a\neq r$) and outgoing arcs
$C_a\setminus{s_a}$, so its block is $C(\nu_a)=C_a$; the added $\ell_x$ are
leaves and contribute no block.  Every maximal clique of $G$ has at least two
vertices ($G$ is connected with $p\ge2$), so $|C(\nu_a)|\ge 2$ for all $a$, and
in particular the root $\nu_r$ has out-degree $|C_r|\ge 2$.  Equipping $T$ with
any $\rho$ satisfying Assumption~\ref{ass:pos} and applying
Theorem~\ref{thm:main} and Proposition~\ref{prop:blocks}, the maximal cliques
of $\Grph(T)$ are exactly the $C_a$, whence
$\Ee(\Grph(T))=\bigcup_a\binom{C_a}{2}=\Ee(G)$ and $\Grph(T)=G$.
\end{proof}

\begin{remark}
\label{rem:realis-lines}
Conditions (i)–(iii), for a connected $G$, characterise exactly the connected
line graphs of trees (Remark~\ref{rem:linegraph}); Theorem~\ref{thm:realis} is
thus the statement that these are precisely the realisable Markov graphs.  The
realising undirected tree $\Tu$ is unique up to isomorphism
(Proposition~\ref{prop:recovery}); the only remaining freedom is the rooting of
Remark~\ref{rem:rooting}, and the root chosen in the construction above (an
arbitrary clique $C_r$) is just one such choice.  For a disconnected $G$
satisfying (i)–(iii) componentwise, the same argument realises $G$ by a
\emph{forest} of accumulation trees.
\end{remark}

The practical consequence is that graphical compatibility with a tree-based
accumulation mechanism is \emph{decidable} from the graph alone, and cheaply so:
chordality and the maximal cliques of a graph are obtained in linear time by
lexicographic breadth-first search \citep{Rose1976}, after which conditions (ii) and (iii) are
read off from the clique–vertex incidence.  This is the basis of
Section~\ref{sec:diagnostics}.

\section{Implications for inference and structure inference}
\label{sec:inference}

The structural properties established above have direct implications for
inference and model selection, and provide a justification for the use of
decomposable graphical models here that goes beyond computational convenience.

\subsection{Likelihood factorisation}

Since $\Grph(T)$ is decomposable with singleton separators, the joint
distribution of $Y$ factorises as in \eqref{eq:factor}: a product of clique
marginals divided by a product of \emph{univariate} separator marginals.  The
likelihood is therefore a product of low-dimensional terms, each clique
connected to the rest of the model through a single binary variable, and
maximum likelihood estimation is available in closed form clique by clique when
the estimate is interior; with sparse clique tables it may fall on the boundary,
where the closed form is the extended (limiting) maximum likelihood estimate.
The exponential factor in the cost of exact inference is governed by
$\max_a|C_a|$, i.e.\ by the branching degree of the event tree
(Corollary~\ref{cor:edges}(3)); the total cost grows in addition linearly with
the number of cliques and of variables.

\subsection{Model dimension: a three-tier hierarchy}
\label{sec:dimension}

It is worth being precise about how much of the accumulation structure the
graphical model actually captures, because the answer bears directly on the
diagnostics of Section~\ref{sec:diagnostics}.

By Remark~\ref{rem:bijection} the accumulation model is a saturated multinomial
on the $p+1$ nodes of $T$, and therefore has exactly
\[
  \dim \mathcal{M}_{\mathrm{acc}} \;=\; p
\]
free parameters. This dimension decomposes along the junction
tree.  In the accumulation model the block $C(v)$ admits only $|C(v)|+1$
configurations of positive probability, either the lineage has not reached
$v$ (then every variable of $C(v)$ is zero; this case does not arise for the
root), or it has reached
$v$ and either stopped there or taken exactly one of the $d(v)$ outgoing
branches, so the corresponding marginal has $|C(v)|$ free parameters.  Using
Corollary~\ref{cor:edges}(2),
\[
  \sum_{a=1}^{m}|C_a| \;-\; (m-1) \;=\; \bigl(p+m-1\bigr)-(m-1) \;=\; p ,
\]
exactly the dimension of $\mathcal{M}_{\mathrm{acc}}$, with the $m-1$ subtracted
terms being the singleton separators.

Contrast this with the graphical models it sits inside.  Writing
$\mathcal{M}_{\mathrm{sing}}$ for the (unrestricted) binary decomposable model
on $\Grph(T)$ with singleton separators, the strictly positive distributions on
$\{0,1\}^p$ that factorise over $\Grph(T)$, and $\mathcal{M}_{\mathrm{sat}}$ for
the saturated model of all strictly positive distributions on $\{0,1\}^p$,
\[
  \dim\mathcal{M}_{\mathrm{acc}} = p,
  \qquad
  \dim\mathcal{M}_{\mathrm{sing}} = \sum_{a=1}^{m}\bigl(2^{|C_a|}-1\bigr)-(m-1),
  \qquad
  \dim\mathcal{M}_{\mathrm{sat}} = 2^{p}-1,
\]
and $\mathcal{M}_{\mathrm{acc}}\subset\overline{\mathcal{M}}_{\mathrm{sing}}
\subset\overline{\mathcal{M}}_{\mathrm{sat}}$, where the bar denotes closure:
the accumulation distributions factorise over $\Grph(T)$ by \eqref{eq:factor},
but they have structural zeros, so they are limits of members of
$\mathcal{M}_{\mathrm{sing}}$ rather than members.  For the tree of Figure~\ref{fig:example}
($p=7$, cliques of sizes $2,4,3$) these are $7$, $23$ and $127$ respectively.

Two consequences follow.  First, \emph{the graphical model is a relaxation}:
Theorem~\ref{thm:main} characterises the conditional independence structure of
the accumulation process, but the process imposes further, non-graphical
restrictions, the deterministic implications $Y_i=1\Rightarrow Y_j=1$ for
ancestral $j$, which restrict the support of the clique marginals rather than
their conditional independences.  Compatibility with $\Grph(T)$ is therefore
\emph{necessary but not sufficient} for the accumulation mechanism, and any
diagnostic based on the graph alone should be described as such.  Second, and
this shapes the diagnostics of Section~\ref{sec:lrt}, because
$\mathcal{M}_{\mathrm{acc}}$ imposes structural zeros it lies on the
\emph{boundary} of $\mathcal{M}_{\mathrm{sing}}$ rather than in its interior.  The
gap in dimension is therefore \emph{not} a regular number of degrees of freedom,
and the two models cannot be compared by an ordinary likelihood-ratio test; the
appropriate comparison is a test of support compatibility, discussed in
Section~\ref{sec:lrt}.

\subsection{Structure inference}

Because the true Markov graph is decomposable with singleton separators,
searches over general undirected graphs may be restricted to this class
\emph{when the accumulation model is taken as given}.  We stress that this is a
modelling assumption and not, as it is sometimes phrased, a restriction
``without loss of generality'': the whole point of Section~\ref{sec:diagnostics}
is that the restriction is refutable, and it can only be refuted by a procedure
that is allowed to leave the class.  In practice, one therefore fits within the
restricted class and compares against a wider one.

Within the class, the search space is highly structured (though it may still be
combinatorially large).  By Theorem~\ref{thm:realis} the admissible graphs are
exactly the connected line graphs of trees, chordal, with pairwise singleton
clique intersections and no vertex in three cliques; by
Corollary~\ref{cor:edges} such a graph is determined, up to relabelling, by its
junction tree together with the clique sizes (on a fixed set of variables one
must in addition assign the separator and private variables).  Forest- and tree-based selection
procedures \citep{Edwards2010, Abreu2010} are natural here, and the treewidth
bound of Corollary~\ref{cor:edges}(3) guarantees that scoring remains tractable
whenever the branching degree is bounded.

\section{Instantiation: mutation accumulation and lineage tracing}
\label{sec:lineage}

Everything above is independent of any application, and a reader who is content
with the statistical formulation of Section~\ref{sec:terminology}, a single
draw from a tree, recorded through the indicator vector of a laminar family, may
stop here.  This section is the first at which subject-matter knowledge is used,
and it exists to answer one question: is there a real setting in which
Assumption~\ref{ass:pos} is \emph{substantively} correct rather than a technical
convenience?  There is, and identifying it is the point of the section.

\subsection{Mapping events to mutations}

Each event $E_i$ corresponds to a mutation at a specific site.  The tree $T$ is
the mutation tree relating the sampled units, the root being a common ancestral
state.  A mutation is assumed to arise at most once on the tree and, once it has
arisen, to be inherited by every descendant of the point at which it arose, so
that the mutations carried by a sampled unit are exactly those on the path from
the root to its position in the tree.  The binary vector $Y=(Y_1,\dots,Y_p)$
recording presence or absence of mutations in a sampled unit is therefore
generated exactly by the process of Section~\ref{sec:model}, and irreversibility
and persistence in the sense of Section~\ref{sec:model} are the biological
assumptions of no reversion and faithful inheritance.

\subsection{Why lineages are observed at arbitrary stages}

The essential point concerns \emph{what is sampled}.  In the classical
phylogenetic picture, sequences are contemporaneous taxa at the tips of a tree,
and each observation is a complete root-to-leaf path.  By
Proposition~\ref{prop:leafcollapse}, that sampling scheme destroys the very
structure we wish to exploit: the internal arc indicators become deterministic
functions of the others and the Markov graph collapses to a disjoint union of
sibling-leaf cliques.  Leaf sampling is thus the \emph{wrong} idealisation for
the present purpose, and it is worth being explicit that this is a statement
about the observation scheme, not about the biology.

Lineage tracing has a different structure.  A \emph{haplotype} is the observed
combination of mutations carried by a sampled unit, in the notation of this
paper, a realisation of the vector $Y$.  The decisive empirical fact is that in
lineage-tracing data \emph{some sampled haplotypes are ancestral to other
sampled haplotypes}: the sample can contain units caught at intermediate stages
of accumulation, not only at the end.  In mitochondrial lineage tracing and in
single-cell mitochondrial genotyping, intermediate states are observed directly
and routinely \citep{Ludwig2019}; and in mtDNA haplotype networks used in human
identification and population genetics, ancestral haplotypes persist in the
sample alongside the derived haplotypes that descend from them
\citep{Amorim2019}.  In such data the observation point of a lineage is
plausibly modelled as a random node of $T$, internal nodes included, as in
Assumption~\ref{ass:pos}, though the assumption of positive mass at \emph{every}
node is stronger than the mere presence of some ancestral haplotypes in the
sample.  A similar structure arises for somatic mutation accumulation in tumour
progression, where sampled cells may span a range of stages of the progression,
and in cultural or linguistic transmission, where conservative and innovative
variants of a tradition coexist in the same sample.

\subsection{Model assumptions}

The assumptions of Section~\ref{sec:model} correspond to standard simplifying
assumptions in molecular evolution.  That each event occurs at most once is the
assumption of \emph{no recurrent mutation}: a site does not mutate twice
independently in different parts of the tree.  Persistence, an event, once
realised, is retained by every descendant, is the absence of \emph{reversion},
that is, of a back-mutation restoring the ancestral state.  The restriction to
a single tree reflects the absence of \emph{recombination}: no unit inherits
material from two distinct ancestral histories, so every unit has a single
line of descent and the histories form a tree rather than a network.  Taken
together, these three assumptions, no recurrent mutation, no reversion and no
recombination, are precisely the \emph{infinite-sites} model of molecular
evolution, under which the data admit a lossless representation as a
\emph{perfect phylogeny} \citep{Gusfield1991}; the model of the present paper is
its single-draw, arbitrary-stage observation counterpart.  Assumption~\ref{ass:pos}
is, however, \emph{stronger} than a bare perfect phylogeny: it requires every
intermediate genotype to carry positive sampling probability, and, via the
one-event-per-arc convention, it resolves several mutations falling between two
sampled states into an ordering with intermediate nodes that need not themselves
be observed.  The present setting is thus an infinite-sites/perfect-phylogeny
model equipped with a full-node-support, arbitrary-stage observation scheme.
These are idealisations, but they are commonly
adopted where mutation rates are low, and recombination is negligible, in
particular for mitochondrial genomes, chloroplast genomes, and many bacterial
genomes \citep{Didelot2010}.  In such settings, the dependence structure among
mutation indicators is governed by the constraints of
Theorem~\ref{thm:main}, and the departures catalogued in
Section~\ref{sec:diagnostics} are typical consequences of failures of these
assumptions, although a given departure need not identify which assumption has
failed.

\section{Diagnostics and departures from the model}
\label{sec:diagnostics}

\subsection{Graphs incompatible with event accumulation}
\label{sec:incompat}

Theorem~\ref{thm:realis} converts the structural theory into an exact
compatibility criterion.  We stress (Section~\ref{sec:dimension}) that
graph-compatibility is \emph{necessary but not sufficient} for the mechanism: the
criterion tests the graphical constraints of Theorem~\ref{thm:main}, not the
accumulation model in full.  An estimated graph $\widehat{\Grph}$ is
\emph{graph-compatible} with irreversible event accumulation on a \emph{single}
tree if and only if it is connected, chordal, its maximal cliques meet pairwise
in at most one vertex, and no vertex lies in three or more maximal cliques.
Define the two \emph{defects}
\[
  \delta_{\cap}(\Grph)
  \;=\;
  \max\Bigl\{\,0,\ \max_{C_a\neq C_b}\bigl|C_a\cap C_b\bigr|-1\,\Bigr\},
  \qquad
  \delta_{\deg}(\Grph)
  \;=\;
  \max\Bigl\{\,0,\ \max_{v\in\Vv}\ \bigl|\{a : v\in C_a\}\bigr|-2\,\Bigr\},
\]
with the convention $\delta_{\cap}=0$ when $\Grph$ has a single maximal clique.
Then $\Grph$ is graph-compatible with a single tree if and only if it is
connected and chordal with $\delta_{\cap}(\Grph)=\delta_{\deg}(\Grph)=0$; a
disconnected graph meeting the three clique conditions is realisable by a
\emph{forest} rather than one tree (Remark~\ref{rem:realis-lines}).

Both defects are needed.  Figure~\ref{fig:incompatible}(a) shows a graph with
$\delta_{\cap}=1$: the clique overlap $\{Y_2,Y_3\}$ has size two, and no
accumulation tree can produce it, because a block intersection is the single
arc joining two nodes.  Figure~\ref{fig:incompatible}(b) shows a graph that is
chordal, has $\delta_{\cap}=0$, and even admits a junction tree with singleton
separators, ,yet has $\delta_{\deg}=1$, since the hinge vertex lies in three
maximal cliques.  It is therefore \emph{also} incompatible, though a diagnostic
based on clique overlaps alone would accept it (cf.\ Remark~\ref{rem:threecliques}).

\begin{figure}[t]
\centering
\begin{tikzpicture}[
   yv/.style={circle,draw,inner sep=1pt,minimum size=6.5mm,font=\small}]

\begin{scope}
\node[yv] (a1) at (1.5,2.6) {$Y_1$};
\node[yv] (a2) at (0.3,1.4) {$Y_2$};
\node[yv] (a3) at (2.7,1.4) {$Y_3$};
\node[yv] (a4) at (0.6,0.0) {$Y_4$};
\node[yv] (a5) at (2.4,0.0) {$Y_5$};
\draw (a1)--(a2); \draw (a1)--(a3); \draw (a2)--(a3);
\draw (a2)--(a4); \draw (a2)--(a5); \draw (a3)--(a4); \draw (a3)--(a5);
\draw (a4)--(a5);
\node at (1.5,-1.0) {\small (a) $\delta_{\cap}=1$: overlap $\{Y_2,Y_3\}$};
\end{scope}

\begin{scope}[xshift=8.0cm]
\node[yv] (x)  at (1.5,1.3) {$Y_1$};
\node[yv] (b1) at (0.0,2.5) {$Y_2$};
\node[yv] (b2) at (0.0,1.1) {$Y_3$};
\node[yv] (b3) at (3.0,2.5) {$Y_4$};
\node[yv] (b4) at (3.0,1.1) {$Y_5$};
\node[yv] (b5) at (1.0,-0.1) {$Y_6$};
\node[yv] (b6) at (2.4,-0.1) {$Y_7$};
\draw (x)--(b1); \draw (x)--(b2); \draw (b1)--(b2);
\draw (x)--(b3); \draw (x)--(b4); \draw (b3)--(b4);
\draw (x)--(b5); \draw (x)--(b6); \draw (b5)--(b6);
\node at (1.5,-1.0) {\small (b) $\delta_{\deg}=1$: $Y_1$ in three cliques};
\end{scope}
\end{tikzpicture}
\caption{Two graphs incompatible with irreversible event accumulation on a
tree.  (a) The clique overlap $\{Y_2,Y_3\}$ has size two, violating
Theorem~\ref{thm:main}(2).  (b) The graph is chordal, all clique intersections
are singletons, and it admits a junction tree with singleton separators; but
the hinge vertex $Y_1$ lies in three maximal cliques, violating
Theorem~\ref{thm:main}(3).  An arc $(u,v)$ belongs to at most two blocks, ,as an
outgoing arc of $u$ and as the incoming arc of $v$, ,so (b) cannot arise from
any accumulation tree.  A diagnostic based on clique overlaps alone would fail
to detect this.}
\label{fig:incompatible}
\end{figure}

\subsection{Likelihood-based diagnostics}
\label{sec:lrt}

The structural constraints suggest a direct statistical diagnostic.  Given
data, one may use likelihood
ratios (or penalised criteria such as AIC or BIC; \citealp{Edwards2010}) to compare the fit of decomposable models constrained to singleton
separators with that of unrestricted decomposable models.  If the
best singleton-separator model fits comparably to the unrestricted one, the
data are consistent with the graphical constraints of tree-based irreversible accumulation; a substantial
improvement under larger separators indicates dependence that the simple
mechanism cannot produce.

Ordinary likelihood-ratio theory applies to the comparison of two \emph{fixed}
nested decomposable models when the true distribution is strictly positive and
lies in the smaller model (an interior point); the degrees of freedom are then
read off from the dimension formula: for a decomposable binary model
with maximal cliques $C_1,\dots,C_m$ and separators $S_2,\dots,S_m$,
\[
  \dim
  \;=\;
  \sum_{a=1}^{m}\bigl(2^{|C_a|}-1\bigr)
  \;-\;
  \sum_{b=2}^{m}\bigl(2^{|S_b|}-1\bigr),
\]
reducing to $\sum_a(2^{|C_a|}-1)-(m-1)$ under singleton separators.

It is tempting to go one step further and test the $p$-dimensional accumulation
model $\mathcal M_{\mathrm{acc}}$ directly against the singleton-separator model
$\mathcal M_{\mathrm{sing}}$, in whose closure it lies, on the nominal
$\dim\mathcal M_{\mathrm{sing}}-p$ degrees of freedom.  This does \emph{not} give
an ordinary likelihood-ratio test.  As noted in Section~\ref{sec:dimension},
$\mathcal M_{\mathrm{acc}}$ has structural zeros and lies on the \emph{boundary}
of $\mathcal M_{\mathrm{sing}}$, so the comparison is non-regular; for clean data
the statistic is in fact degenerate.  If every observed pattern is an admissible
root-to-node configuration, the empirical distribution itself lies in the
closure of $\mathcal M_{\mathrm{acc}}$ (it is a distribution on the nodes,
possibly with some node masses zero, rather than an interior member of
$\mathcal M_{\mathrm{acc}}$) and factorises over $\Grph(T)$, so the two maximised
likelihoods coincide and the statistic is \emph{zero}; if even one forbidden
pattern is observed, the clean accumulation model assigns it likelihood zero and
the statistic is \emph{infinite}.  The appropriate clean-data diagnostic is
therefore not a likelihood ratio but a \emph{support-compatibility} test: does
every observed binary vector correspond to a root-to-node path in a common event
tree?  Under observation error, one embeds the process in a positive-probability
misclassification model and compares within that latent-variable family, whose
asymptotics must be derived for the model at hand; the nominal
degrees-of-freedom count does not transfer.

In practice, the branching skeleton is unknown and must be estimated; a workable
recipe is to fit a decomposable or forest graph by a penalised search, the
minimal-AIC/BIC search of \citet{Abreu2010, Edwards2010} ,and then test, as
above, whether its separators may be taken to be singletons and whether any
vertex lies in three cliques (Theorem~\ref{thm:realis}).  Estimation on the
resulting decomposable models is governed by the hyper Markov theory of
\citet{DawidLauritzen1993}; when observation noise is present, the
hidden-variable elaboration of Section~\ref{sec:noise}, a latent accumulation
process observed through an error-prone channel, as in the hidden conjunctive
Bayesian network of \citet{Gerstung2009}, is the natural vehicle.

\subsection{A worked example}
\label{sec:worked}

To exhibit the diagnostic in action we simulate lineages from the tree of
Figure~\ref{fig:example} under Assumption~\ref{ass:pos}, drawing the observation
node $Z$ from a strictly positive distribution over $\Nodes$.  As a departure we
introduce a single \emph{recurrent event}: on lineages that follow the $E_6$
branch, event $E_5$ is allowed to occur a second time, independently of the arc
to which it is nominally attached, exactly the parallel or recurrent evolution
catalogued above.  It couples $Y_5$ to $Y_6$, and hence to $Y_1$, without passing
through the hinge $Y_4$.

At the population level the effect is unambiguous, and instructive.  The exact
Markov graph of the tree mechanism is the chordal, singleton-separator graph of
Figure~\ref{fig:example}(b); the recurrent-event mechanism adds the edge
$Y_5\!-\!Y_6$ and thereby creates the chordless four-cycle
$Y_6\!-\!Y_5\!-\!Y_4\!-\!Y_1\!-\!Y_6$.  Its maximal cliques are
$\{Y_1,Y_2,Y_3,Y_4\}$, $\{Y_4,Y_5,Y_7\}$, $\{Y_1,Y_6\}$ and $\{Y_5,Y_6\}$, which
still meet pairwise in at most one vertex, and no vertex lies in three of them, so
$\delta_{\cap}=\delta_{\deg}=0$.  The graph is nonetheless \emph{not chordal}, and
the departure is detected by the chordality condition alone, a concrete
illustration that clause~(1) of Theorem~\ref{thm:main} is not implied by the
clique conditions (ii)--(iii) and must be checked in its own right.

In finite samples we use a BIC edge-addition diagnostic, not a classical score
test, anchored at the hypothesised tree graph $\Grph(T)$.  For each candidate
edge $e$ whose addition keeps the graph chordal, let $\mathrm{BIC}_e$ be the BIC
of the fitted decomposable model and $\mathrm{BIC}_0$ that of $\Grph(T)$; the
statistic is
\[
  D \;=\; \max\Bigl\{\,0,\ \max_{e}\bigl(\mathrm{BIC}_0-\mathrm{BIC}_e\bigr)\Bigr\},
\]
and $D>0$ flags a departure.  Note that the true population edge $Y_5\!-\!Y_6$ is
itself \emph{not} chordality-preserving, so the search cannot add it; it selects
instead a chord (here almost always $Y_1\!-\!Y_5$) that best improves the
decomposable approximation to the non-chordal truth.  The diagnostic is thus a
proxy for the chordality violation rather than a direct test of the offending
edge.

\begin{figure}[t]
\centering
\includegraphics[width=0.66\textwidth]{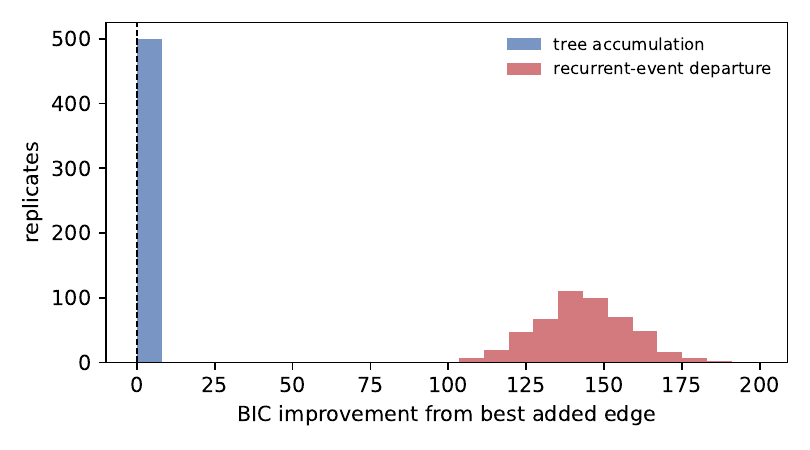}
\caption{Finite-sample behaviour of the singleton-separator diagnostic
(Section~\ref{sec:worked}), over $500$ simulated data sets of size $n=1000$ from
the tree of Figure~\ref{fig:example}.  For each data set, the statistic is the
largest BIC decrease obtainable by adding one chordality-preserving edge to the
hypothesised tree graph; a positive value flags a departure.  Under tree
accumulation the statistic is $0$ in every replicate; under a recurrent-event
departure it is large in every replicate.}
\label{fig:sim}
\end{figure}

We draw $Z$ from the node distribution proportional to
$(13,11,7,9,17,5,19,23)$ on $(z_1,\dots,z_8)$, set the recurrence probability to
$r=0.6$, and take $n=1000$ observations; the BIC uses the dimension formula
above.  Over $500$ replicates,
$D>0$ occurred in $0/500$ replicates under the tree mechanism and in $500/500$
under the recurrent-event mechanism, with median $D\approx 142$;
Figure~\ref{fig:sim} shows the complete separation of the two empirical
distributions.  Consistently with Section~\ref{sec:noise}, independent
misclassification noise applied to the tree mechanism is \emph{also} flagged: at
$2\%$ per-coordinate noise $D>0$ in about a quarter of replicates, because a noisy
hinge no longer screens exactly, and the observed graph fills in.
The clean diagnostic thus
tests accumulation \emph{and} clean observation jointly.

\subsection{Recombination and blockwise models}

The results rely on the assumption that the observed variables are generated by
accumulation along a \emph{single} tree.  Recombination and horizontal transfer
violate this by merging two distinct ancestral histories into one unit, so
that the histories form a network rather than a tree \citep{Didelot2010}.  One way to accommodate
them is through blockwise models, in which different subsets of variables are
associated with different underlying trees, with dependence between blocks
mediated by a limited number of shared variables.  Within each block the
conclusions of Theorem~\ref{thm:main} apply.  We emphasise, however, that the
\emph{global} graph in such a model need not be decomposable at all, and that a
mixture of decomposable models is in general not a graphical model; the correct
statement is that the block-and-separator structure is inherited locally, and
that the scale at which decomposability breaks down is itself informative about
the scale of recombination.  A systematic treatment is beyond the scope of this
paper.

\subsection{Observation noise}
\label{sec:noise}

Suppose the observed vector is $\tilde Y$, obtained from the latent $Y$ by
independent misclassification, $\Pr(\tilde Y_i\neq Y_i)=\varepsilon_i\in(0,1)$,
independently across $i$ given $Y$.  It is tempting to assume that noise
leaves the conditional independence structure intact and merely weakens the
observed dependencies.  For the \emph{latent} vector this is trivially true;
for the \emph{observed} vector it is false, and in the wrong direction.

The screening mechanism of Lemma~\ref{lem:hinge} operates through $Y_h$
\emph{exactly}: conditioning on the true hinge value renders one of the two
blocks degenerate.  Conditioning on a noisy surrogate $\tilde Y_h$ does not,
since $\tilde Y_h$ leaves residual uncertainty about $Y_h$, through which
dependence continues to propagate.  Consequently, the Markov graph of $\tilde Y$
is typically \emph{denser} than that of $Y$, and the singleton-separator
structure is the first casualty, not the last.  For the tree of
Figure~\ref{fig:example}, for instance, the Markov graph of $\tilde Y$ is
complete for every common misclassification rate we tried (from $0.02$ to
$0.1$) and for unequal rates, and noise on the hinge $Y_4$ alone already adds
six edges.  We do not prove that completeness holds generically; what the
argument above establishes is that exact hinge screening is in general
destroyed.  What survives is a matter of
magnitude, not of structure: for small $\varepsilon_i$ the spurious edges carry
weak partial associations, and the block-and-separator pattern remains
detectable as the dominant signal, but it is no longer exact.

The practical implication is the opposite of the reassuring one.  A graph
estimated from noisy accumulation data will tend to violate
Theorem~\ref{thm:main} even when the mechanism is correct, so the diagnostics of
Section~\ref{sec:incompat} must be applied with an explicit noise model, or
interpreted as tests of ``accumulation \emph{and} clean observation'' jointly.
Latent-variable treatments that model $Y$ explicitly and integrate out the
misclassification recover the exact structure (for example the hidden conjunctive
Bayesian network of \citealp{Gerstung2009}), at the price of an EM or
Bayesian computation; we regard this as the natural next step and do not pursue
it here.

\subsection{Approximate accumulation}

The model assumes strict irreversibility and persistence.  In practice these may
hold only approximately: events may occasionally be lost, reverted, or
misrecorded.  Such departures can be viewed as perturbations of the structural
assumptions, under which the graphical properties established here describe the
dominant dependence structure, with deviations from irreversibility introducing
additional and typically weaker interactions.  When departures are rare, the strong dependencies are still organised into
blocks connected by hinge variables.  We caution, however, that the exact
conditional-independence graph is \emph{not} a continuous function of the
distribution: an arbitrarily small generic perturbation can add every missing
edge, so the induced graph need not stay close to a singleton-separator graph.
What remains close is the \emph{magnitude} pattern, the partial associations or
conditional mutual information, rather than the exact graph.  A formal
theory, in particular, rates at
which reversion or event loss degrades the graphical structure, would be a
natural companion to the noise analysis of Section~\ref{sec:noise}.

\section{Concluding remarks}
\label{sec:discussion}

This paper establishes a direct connection between irreversible
event-accumulation processes on trees and decomposable graphical models.  The
conditional independence structure induced by such processes is highly
constrained: the Markov graph is chordal, its maximal cliques correspond to the
nodes of the underlying undirected event tree of degree at least two (so that
$\Grph(T)=\Lg(\Tu)$ is the line graph of that tree), no vertex lies in three
cliques, and the separators are the singleton arc-events joining two such nodes.
A single screening identity, each hinge variable separates its subtree from the
rest, delivers the entire structure.  The converse (Theorem~\ref{thm:realis})
shows that these constraints are not merely necessary but characterise the
realisable graphs exactly, so that graphical compatibility with the mechanism
is decidable from the estimated graph alone.

Three lessons stand out.  First, \emph{the modelling of the observation scheme
matters}: lineages must be observed at arbitrary stages of progression, not only
at their endpoints, or the flow constraints of complete lineages
collapse the graph entirely (Proposition~\ref{prop:leafcollapse}).  Second,
chordal independence graphs, often introduced for mathematical convenience
(because they yield decomposable models with tractable likelihoods), also arise
naturally as representations of irreversible event accumulation.  Third, the
graphical model is a \emph{relaxation} of the accumulation model, of much larger
dimension (Section~\ref{sec:dimension}); the graph captures the conditional
independences but not the support restrictions that irreversibility imposes.
When an inferred graph violates the constraints derived here, the deviation is
itself informative. It signals recombination, recurrent events, mixtures of
distinct populations, or observational noise.  More broadly, the results illustrate how
structural properties of graphical models can reflect fundamental constraints
imposed by the underlying generative mechanism in any system governed by the
irreversible accumulation of events along lineages.

\end{document}